\documentclass[12pt]{article}
\usepackage{graphicx}
\usepackage{mathptmx}
\usepackage{latexsym,amsmath,amssymb,amsfonts,amsthm}

\newcommand{\ti}[1]{\mbox{\tiny $#1$}}
\newcommand{\scr}[1]{\mbox{\scriptsize $#1$}}
\newcommand{\ft}[1]{\mbox{\footnotesize $#1$}}
\newcommand{\sm}[1]{\mbox{\small $#1$}}
\newcommand{\la}[1]{\mbox{\large $#1$}}
\newcommand{\La}[1]{\mbox{\Large $#1$}}

\newfont{\lie}{eufm10 at 12pt}
\newfont{\field}{msbm10 at 11pt}

\newtheorem{theorem}{Theorem}[section]
\newtheorem{lemma}{Lemma}[section]
\newtheorem{corollary}{Corollary}[section]
\newtheorem{proposition}{Proposition}[section]

\begin{document}
\title{
Cauchy-Riemann inequalities on $2$-spheres of $\mathbb{R}^7$}
\markright{\sl 
\hfill Cauchy-Riemann inequalities \hfill}
\author{ Isabel M.C.\ Salavessa}
\date{}
\protect\footnotetext{\!\!\!\!\!\!\!\!\!\!\!\!\! {\bf MSC 2000:}
Primary:  53C42; 53C38. Secondary: 58E12; 35J19; 47A75
{\bf ~~Key Words:} Stability,  Parallel
Mean curvature,  Calibration,  Cauchy-Riemann inequality. \\
Partially supported by
FCT through programs
 PTDC/MAT/101007/2008, PTDC/MAT/118682/2010.}
\maketitle ~~~\\[-10mm]
{\footnotesize Centro de F\'{\i}sica das Interac\c{c}\~{o}es
Fundamentais, Instituto Superior T\'{e}cnico, Technical University
of Lisbon, Edif\'{\i}cio Ci\^{e}ncia, Piso 3, Av.\ Rovisco Pais,
1049-001 Lisboa, Portugal;~ isabel.salavessa@ist.utl.pt}\\[5mm]
{\small {\bf Abstract:} We prove that an integral Cauchy-Riemann inequality
holds for any pair of smooth functions $(f,h)$ on the 2-sphere
 $\mathbb{S}^2$, and equality
holds iff $f$ and $h$ are related $\lambda_1$-eigenfunctions.
We extend such  inequality to  $4$-tuples of functions, 
only valid  on the $L^2$-orthogonal complement of 
a suitable nonzero finite dimensional space 
of functions. As a consequence we prove that  
$2$-spheres are not $\Omega$-stable surfaces with parallel mean curvature
in $\mathbb{R}^7$ for the associative
calibration $\Omega$. 
}
\section{Introduction}
In \cite{Sa1} we extended to submanifolds with higher codimension
the variational characterization of hypersurfaces
of Riemannian manifolds with constant mean curvature $H$ discovered 
by Barbosa, do Carmo and Eschenburg
\cite{BdC,BdCE}.
Given an $m$-dimensional oriented immersed 
submanifold $\phi:M\to \bar{M}$
  of an $(m+n)$-dimensional calibrated 
Riemannian manifold $(\bar{M}, \bar{g})$ 
with a  semicalibration $\Omega$ of rank $m+1$ (see \cite{PumRiv}, in 
\cite{Sa1} we denominated it by precalibration), 
and  assuming
 $\phi$ has calibrated extended tangent space on a domain $D$, 
that is, along $D$ the vector bundle
$EM=\mathbb{R}\nu \oplus TM$ is  $\Omega$-calibrated,
where $\nu$ is a globally defined unit normal on $D$ such that
$H\in \mathbb{R}\nu$, then we proved that $\phi$ has constant mean
curvature on $D$ if and only if $\phi$ is a critical point
of the area $A_D(\phi)$ for all variations $\bar{\phi}:[0,\epsilon)\times D\to
\bar{M}$  of $\phi$, $\bar{\phi}(t,p)=\phi_t(p)$, $\phi_0=\phi$,
 that fixes the  boundary of $D$ (in case this one exists) and
preserves the $\Omega$-volume
$V_D(t)=\int_{[0,t]\times D}
\bar{\phi}^*\Omega$.
The later condition means  that  $V_D(t)$ is constant on $t$, that is
 $V_D(t)=V_D(0)=0$.
This turns out to be equivalent to $\phi$ to be a critical
point of $J_D(t)=A_D(\phi_t)+mh_DV_D(t)$,
for any variation fixing the boundary of $D$,
 where $h_D$ is the mean value
of $\|H\|$ on $D$.
The second variation of $J_D(t)$ was computed, obtaining
${J''}_D(0)=\int_D\bar{g}(\mathcal{J}_{\Omega}(W), W)dM=:I_{\Omega}(W,W)$, 
where  $W=\frac{d}{dt}_{|_{t=0}}\phi_t$ is the vector variation field
for variations that preserve the $\Omega$-volume, and 
\begin{equation}
\mathcal{J}_{\Omega}(W)=-\Delta^{\bot}W^{\bot}-\bar{R}(W^{\bot})-
\tilde{B}(W^{\bot})+ m\|H\|C_{\Omega}(W^{\bot}).
\end{equation}
This second order differential operator depends only
on the normal component of $W$ and it 
is the usual Jacobi operator with an extra term
$m\|H\|C_{\Omega}(W^{\bot})$, a  $L^2$-self-adjoint first order differential
operator defined for $W\in C_0^{\infty}(NM_{/D}),$ such that
\begin{equation}
\int_D \bar{g}(C_{\Omega}(W),W)=
\int_D \La{(}\sum_i\Omega(W, e_1, \ldots,
\nabla_{e_i}^{\bot}W_{\ti{(i)}}, \ldots, e_m)
+(\bar{\nabla}_{W}\Omega)({W},e_1, \ldots, e_m) \La{)}dM,
\end{equation}
where $e_1,\ldots,e_m$ defines a direct o.n. frame of $TM$.
We also called by the same name the operator
$\mathcal{J}'_{\Omega,D}(W)=\mathcal{J}_{\Omega}(W)
-\Psi_{\Omega, D}(W)\nu$, where $\Psi_{\Omega, D}(W)$ is the linear
operator:
$$\Psi_{\Omega, D}(W)=\frac{1}{|D|}\int_{D}\bar{g}(\mathcal{J}_{\Omega}(W),\nu)
dM$$
For simplicity we assumed $\phi$ to have parallel mean curvature, 
and so $\nu$ is a parallel unit normal.
We can extend the definition of $I_{\Omega}(W,W)$ to
the space $H^1_{0,T}(NM_{/D})$ given by
the $H^1$-completion of the vector space
generated by the set of  normal vector  variations $W$
of $\Omega$-volume preserving variations $\bar{\phi}$ fixing the
boundary of $D$, or equivalently, it is the subspace
of $H^1_0(NM_{/D})$ of the normal sections  that satisfy a zero
 mean value property
$$ \int_D\Omega(W,e_1, \ldots, e_m)dM=0.$$
Then we say that $\phi$ is $\Omega$-stable on $D$ if
$I_{\Omega}(W,W)\geq 0$, for all $W\in H^1_{0,T}(NM_{/D})$.
The space $H^1_{0,T}(NM_{/D})$ is just $H^1_{0,T}(D)\oplus H^1_0(F)$
where $H^1_{0,T}(D)=H^1_0(D)\cap L^2_T(D)$ with
$L^2_T(D)=\{f\in L^2(D): \int_DfdM=0\}$, and $F$ is the
normal  subbundle orthogonal complement of $\nu$.
A Morse index theorem can be stated for submanifolds with parallel
 mean curvature and calibrated extended tangent 
space (see Remark 4.3 of \cite{Sa1}).
If $\bar{M}=\mathbb{R}^{m+n}$ and $M$ is closed,
and supposing that $\bar{\nabla}\Omega=0$ or
$C_{\Omega}=0$ (in fact it is sufficient to assume the vanishing of
$\bar{g}(C_{\Omega}(\nu), \nu)$, or equivalently of $\bar{\nabla}_{\nu}
\Omega(\nu, e_1, \ldots, e_m)$), we proved in \cite{Sa1} 
(Theorem 4.2) that under the natural
integral inequality condition $\int_MS(2+h\|H\|)dM\leq 0$,
where $h$ and $S$ are the height functions
$h=\bar{g}(\phi,\nu)$ and $S=\sum_{ij}\bar{g}(\phi, B^F(e_i,e_j))
B^{\nu}(e_i,e_j)$ ($B^{\nu}$ and $B^F$ stand for
  the $\nu$ and $F$-components
of the second fundamental form $B$, respectively), if $M$ is $\Omega$-stable
then $\phi$ is pseudo-umbilical, and in case $NM$ is a trivial bundle,
then $M$ must be a sphere.
So it is a fundamental question to describe for which
semicalibrations $\Omega$ are spheres $\Omega$-stable. If $n=1$
this is completely determined, for in \cite{BdC} it is proved that stable
closed immersed 
hypersurfaces of constant mean curvature are exactly the spheres.
If $n\geq 2$ and $C_{\Omega}=0$ then $m$-spheres $\mathbb{S}^m$
of $\Omega$-calibrated
vector subspaces are $\Omega$-stable. This is also the case
$n=2$ and $\Omega$ parallel (see Corollary 2.1), or for any $n\geq 2$ and  
$\Omega$ is
a semicalibration defined by a   fibration of $\mathbb{R}^{n+m}$
with an $(m+1)$-dimensional totally geodesic fibre where 
$\mathbb{S}^m$ lies (\cite{Sa1}). On the other hand $C_{\Omega}$ does not vanish
for most well known calibrations, namely the ones coming from special 
holonomy. 
In this paper we consider one of such case of a parallel
calibration with non vanishing $C_{\Omega}$, namely the
\em associative calibration, \em   defined  by the $G_2$ structure
of the  Euclidean space $\mathbb{R}^7$, and that is given by the 3-form
$$\Omega= dx^{123}+ dx^{145}+dx^{167}+dx^{246}
-dx^{257}-dx^{347}-dx^{356},$$
and prove that 2-spheres of associative subspaces are $\Omega$-unstable
on  $\mathbb{R}^7$.
The $\Omega$-stability condition is equivalent to the
\em long integral Cauchy Riemann inequality \em  to hold for any $4$-tuples
of functions $f=(f_4,f_5,f_6,f_7):\mathbb{S}^2\to \mathbb{R}^4$:
\begin{eqnarray}
&&4\int_{\mathbb{S}^2}\phi_1 (\langle J\nabla f_4,\nabla f_5\rangle
+\langle J\nabla f_6,\nabla f_7\rangle )dM \nonumber\\
&&+ 4\int_{\mathbb{S}^2}\phi_2 (\langle J\nabla f_4,\nabla f_6\rangle
-\langle J\nabla f_5,\nabla f_7\rangle )dM \nonumber \\
&&- 4\int_{\mathbb{S}^2}\phi_3 (\langle J\nabla f_4,\nabla f_7\rangle
+\langle J\nabla f_5,\nabla f_6\rangle )dM \nonumber \\
&\leq& \int_{\mathbb{S}^2}(\|\nabla f_4\|^2+ \|\nabla f_5\|^2+
\|\nabla f_6\|^2+ \|\nabla f_7\|^2)dM, 
\end{eqnarray}
where $J$ is the complex structure of $\mathbb{S}^2$ and
$\phi=(\phi_1,\phi_2,\phi_3): \mathbb{S}^2\to \mathbb{R}^3$ is the inclusion
map. 
In case two of the functions $f_{\alpha}$ are zero, the  above inequality  
  gives 
 the \em short integral Cauchy Riemann inequalities \em holding
for any pair of functions $(f,h):\mathbb{S}^2\to \mathbb{R}^2$
\begin{equation}
2 \int_{\mathbb{S}^2} \phi_i\langle J\nabla f, 
\nabla h\rangle dM
\leq \sqrt{\int_{\mathbb{S}^2}\|\nabla f\|^2dM }\, 
\sqrt{\int_{\mathbb{S}^2}\|\nabla h\|^2dM} \mbox{~~~~~for~}i=1,2,3
\end{equation}
\noindent
as we can easily see from (3) by replacing $f$ and $h$ by $tf$ and $t^{-1}h$,
respectively, 
where $t^2=\|\nabla h\|_{L^2}/\|\nabla f\|_{L^2}$.
We state our  main results in Theorems 1.1 and 1.2:
\begin{theorem} For each $i\in \{1,2,3\}$ 
the short integral Cauchy Riemann inequality (4) holds for any
smooth maps $f,h: \mathbb{S}^2\to \mathbb{R}$. Furthermore,   
equality  holds in (4) if
and only if either $f$ or $h$ is constant, or $ f= c_j \phi_j+c_k\phi_k +c$
and $h= -c_k\phi_j+ c_j\phi_k+ c'$, 
where $c_j,c_k, c, c'$ are constants,
and $(i,j,k)$ is a positive  permutation of $(1,2,3)$.
\end{theorem} 
\noindent
The eigenvalues for the closed Dirichlet problem on the unit
2-sphere $\mathbb{S}^2$
constitute an increasing sequence $0=\lambda_0<\lambda_1<\ldots< \lambda_l
< \cdots$ converging to infinity. 
We denote by $E_{\lambda_l}$ the eigenspace of dimension $2l +1$
corresponding to the 
 eigenvalue  $\lambda_l=l(l+1)$
and by $E_{\lambda_l}^+$
the union of all eigenspaces $E_{\lambda}$ with $\lambda\geq \lambda_l$.

\begin{theorem}
The long integral Cauchy Riemann inequality holds for all
$f_{\alpha}\in E_{\lambda_0}\oplus 
E_{\lambda_1}\oplus E_{\lambda_2}\oplus E^+_{\lambda_6}$.
But there exist a 4-tuple of functions $f_{\alpha}\in E_{\lambda_3}$
for which (3) is not satisfied. In particular
$\mathbb{S}^2$ is $\Omega$-unstable.
\end{theorem}
\noindent
In Proposition 3.2 we give one more class of functions for which (3) 
is satisfied, as an immediate consequence of Theorem 1.1.
The index is the dimension of the largest vector space of 4-tuples
of functions for which  inequality (3) does not hold.
The exact index will be computed somewhere later. The search of directions
of instability  defined by  4-tuples 
requires long computations. 
Since the functions $f_{\alpha}$ can be expressed as  a $H^1$-sum  of    
 spherical harmonics on $\mathbb{S}^2$, we 
express $f_{\alpha}$ in a non-orthonormal sum of  monomial functions. 
We  show that the  use of eigenfunctions reduces the study of the 
Cauchy-Riemann inequalities
to consider functions $f_{\alpha}$ only in one ore two different eigenspaces.
We prove this by observing
first  that $\int_{\mathbb{S}^2}\phi\langle J\nabla f, \nabla h\rangle$
is a skew-symmetric functional in the three variable functions $(\phi,f,h)$, 
and derive a
Weitzenb\"{o}ck-type formula that concludes that $\langle J\nabla \phi_i,
\nabla f\rangle $ maps a $\lambda_l$-eigenfunction $f$ into
a $\lambda_l$-eigenfunction.
Using only algebraic methods we determine
 that (3) holds 
  for functions in $E_{\lambda_l}$ where $l=0,1,2$ or $l\geq 6$,
but for $l=3$ we need to use   Mathematica and Fortran programming
to diagonalize   a $40\times 40$ matrix, 
to obtain all stable and unstable directions. 
The cases $l=4$, and $l=5$ are considerably more
complicate for it correspond to diagonalize 
a $60\times 60$ and a $80\times 80$ matrix respectively. We do not consider
these two cases here.

Recall that $H^1_0(\mathbb{R}^m)=H^1(\mathbb{R}^m)$. Using the stereographic
projection $\sigma: \mathbb{R}^2\to \mathbb{S}^2$
$$\sigma(w)=\left(\frac{|w|^2-1}{|w|^2+1}, \frac{2 w_1}{|w|^2+1},
\frac{2 w_2}{|w|^2+1}\right),$$
 Theorem 1.1 is translated 
into next Corollary:
\begin{corollary} If $f, h\in H^1(\mathbb{R}^2)$, then for $i=1,2,3$
$$
2\int_{\mathbb{R}^2}\sigma_i
\langle J_0\nabla^0f,
\nabla^0h\rangle dw \leq \|\nabla^0f\|_{L^2}\|\nabla^0h\|_{L^2}$$
where  $J_0$ is the canonical complex structure of 
$\mathbb{R}^2$, and $\nabla^0$ is the gradient operator in $\mathbb{R}^2$.
Furthermore, equality  holds if and only if $f$ or $h$ vanish.
\end{corollary}
\noindent
Note that $\nabla^0h=J_0\nabla^0f$ if and only if $f+ih:\mathbb{R}^2\to
\mathbb{C}$ is an holomorphic map. In this case $f$ and $h$ do not lie
in $H^1(\mathbb {R}^2)$ unless they are zero functions. Furthermore
non-constant holomorphic maps cannot be constant in any open sets, and in
particular on a set
where, for some $i$, $\sigma_i <1-\delta$
with $\delta>0$ small, and so   the coefficient $2$ in the above inequality
( or in (4)) is expectable. Moreover, since $\sigma_i$ are not $L^2$-functions, 
equality only holds if $f$ or $h$ vanish.
\section{Preliminaries }
Let $\bar{M}$ be an $(m+n)$-dimensional Riemannian manifold with a
semicalibration $\Omega$ of rank $(m+1)$. This means that $\Omega$ is
an $(m+1)$-form such that
$ |\Omega(u_1, \ldots, u_{m+1})|\leq 1$
for any o.n. system of vectors $u_i$. We consider
$\phi:M\to \bar{M}$ 
an $m$-dimensional ($m\geq 2$), 
oriented, immersed submanifold with nonzero parallel 
 mean curvature $H=\|H\|\nu$, and calibrated extended
tangent space $H\oplus TM$, that is $\Omega(\nu,e_1, \ldots, e_{m+1})=1$
holds for a d.o.n.\ frame $e_i$ of $TM$.

Given a smooth normal section  $W$ with compact support on a
domain $D$  of $M$, 
 then $C_{\Omega}(W)$ is defined as the normal section such that
for all $W'\in C_0(NM_{/D})$ (cf. \cite{Sa1})
\begin{eqnarray}
\lefteqn{\bar{g}(C_{\Omega}(W),W') =\sum_i\Omega(W', e_1, \ldots,
\nabla_{e_i}^{\bot}W_{\ti{(i)}}, \ldots, e_m)}\\
&&+
\frac{1}{2}\La{(}(\bar{\nabla}_{W}\Omega)({W'},e_1, \ldots, e_m)
+(\bar{\nabla}_{{W'}}\Omega)({W},e_1, \ldots, e_m)\La{)} \nonumber \\
&&-\sum_i \frac{1}{2}\bar{\nabla}_{e_i}
\Omega(W, e_1, \ldots,W'_{\ti{(i)}},\ldots, e_m)
-\sum_{i}\sum_{j\neq i}\frac{1}{2}\Omega(W, e_1, \ldots, 
B(e_i,e_j)_{\ti{(j)}}, \ldots,W'_{\ti{(i)}}, \ldots, e_m)\nonumber
\end{eqnarray}
where  $\scr{(i)}$ means the $i$-
position. A simple computation shows that
\begin{equation}
\bar{g}(C_{\Omega}(W'),W)-\bar{g}(C_{\Omega}(W),W')=div_M(X_{WW'})=
-\delta (\xi(W,W'))
\end{equation}
where $X_{WW'}\in C_0^{\infty}(TM)$ is a vector field on $M$
 and $\xi:\wedge^2NM\to T^*M$  a tensor defined  by 
$$g(X_{WW'}, e_i)=\Omega(W, e_1, \ldots, W'_{\ti{(i)}}, \ldots, e_m)
=\xi(W,W')(e_i)$$
$$ \xi(W,W')(u)=\Omega(W,W',*u),$$
 where $*$ is the star operator on $M$.
Eq.(6) is derived by
taking into consideration that  $M$ has calibrated extended space, and so
$$
{\sum_i\Omega((\bar{\nabla}_{e_i}W)^{\top},
 e_1,\ldots, W'_{\ti{(i)}},\ldots, e_m) =
m\|H\|\bar{g}(W,\nu)\bar{g}(W',\nu)}
$$
which is symmetric on $W,W'$.
Thus, $C_{\Omega}$ is $L^2$-self-adjoint and (2) holds. 
For $n\geq 2$, we recall Lemma 4.4 of \cite{Sa1}. 
We will give here a clearer proof.
\begin{lemma} If $n\geq 2$, 
$C_{\Omega}$ vanish iff (7) and (8) holds:
\begin{eqnarray}
&&\xi \mbox{~vanish~}\\
&&\bar{\nabla}_{W}\Omega(W', e_1, \ldots, e_m)=
-\bar{\nabla}_{W}\Omega(W', e_1, \ldots, e_m)
\end{eqnarray}
If $n=2$ then (7) holds.
\end{lemma}
\noindent
\em Proof. \em  If $C_{\Omega}$ vanish, then both
$\bar{g}$-anti-self-adjoint and  $\bar{g}$-self-adjoint parts 
of $\bar{g}(C_{\Omega}(W),W')$  vanish,
what means $div(X_{WW'})=0$, and by taking normal sections that at a point
have zero normal covariant derivative we see the above equality (8)
 holds.  Then we consider  $\tilde{W}=W+ fW'$ where $\nabla^{\bot}W
=\nabla^{\bot}W'=0$ at a given point $p$ and $f$ is any function. From 
$\bar{g}(C_{\Omega}(\tilde{W}),\tilde{W})=0$  and (5) we conclude that
$\xi(W,W')(\nabla f)=0$ at $p$. Since $f$ is arbitrary we get
$\xi=0$. \qed\\[3mm]

From now on we are assuming $M$ is a $m$-dimensional
Euclidean  sphere  $\mathbb{S}^m_r$ of radius $r$  of a
$\Omega$-calibrated Euclidean subspace $\mathbb{R}^{m+1}$
of  $\mathbb{R}^{m+n}$, and $\phi:\mathbb{S}^m_r\to \mathbb{R}^{m+1}\subset
\mathbb{R}^{m+n}$ denotes the inclusion map, and $\varepsilon_i$, 
$i=1,\ldots,m+1$,  is the canonical basis of $\mathbb{R}^{m+1}$.

We recall that the eigenvalues of $\mathbb{S}^m_r$ for the closed Dirichlet 
problem are given by $\lambda_l(r)=\frac{l(l+m-1)}{r^2}$, 
with $l=0,1, \ldots$, and    the $\lambda_l(r)$-eigenfunctions are
of the form $f_r(x)=f(\frac{x}{r})$ where $f$ is a
$\lambda_l(1)$-eigenfunction of the unit sphere $\mathbb{S}^m$.
We omit the parameter $r$ if $r=1$.
Furthermore, 
if $f\in E_{\lambda_l(r)},~
h\in E_{\lambda_s(r)}$ then $$\int_{\mathbb{S}^m_r}fh\, dM=0 ~~\mbox{if}~
l\neq s \mbox{~~~~and~~~~} 
\int_{\mathbb{S}_r^m}\langle \nabla f, \nabla h\rangle\,  dM= - \delta_{ls}
\lambda_l(r)\int_{\mathbb{S}^m_r}fh\, dM.$$
There exists a  $L^2$-orthonormal basis $\psi_{l,\sigma}$ 
 of $L^2(\mathbb{S}^m_r)$ 
of eigenfunctions
($1\leq \sigma\leq m_l$, where $m_l$ denotes the multiplicity of 
$\lambda_l(r)$).
The Rayleigh characterization of $\lambda_l(r)$ is given by
$$\lambda_l(r)=\inf_{f\in E_{\lambda_l(r)}^{+}}\frac{\int_{\mathbb{S}^m_r}
\|\nabla f\|^2dM}{\int_{\mathbb{S}^m_r}f^2dM},$$ where
$E_{\lambda_l(r)}^{+}$
is the $L^2$-orthogonal complement of the sum of the eigenspaces
$E_{\lambda_i(r)}$, $i=1,\ldots,l-1$. Equality holds for 
$f\in E_{\lambda_l(r)}$. 
Each eigenspace $E_{\lambda_l(r)}$ is exactly composed
 by the restriction to
$\mathbb{S}^m_r$ of the harmonic 
homogeneous polynomials functions of degree $l$ of $\mathbb{R}^{m+1}$, and
it  has dimension $m_l=\binom{m+l}{m}-\binom{m+l-2}{m}$.
Thus, each eigenfunction
$\psi\in E_{\lambda_l(r)}$, is of the form
$\psi=\sum_{|a|=l} \mu_{a }
 \phi^{a}$, where $\mu_{a}$ are
some scalars and 
$a=(a_1,\ldots, a_{m+1})$ denotes a multi-index
of lenght
$|a|=a_1+ \ldots +a_{m+1}=l$ and  
$$\phi^{a}=
\phi_1^{a_1}\cdot\ldots\cdot\phi_{m+1}^{a_{m+1}}.$$
From $\nabla \phi_i= \epsilon_i^{\top}$ and that $\sum_i \phi_i^2=r^2$
we see that
\begin{equation}
\left\{ \begin{array}{l}
\langle \nabla \phi_i,\nabla\phi_j\rangle =\delta_{ij}-\frac{1}{r^2}
\phi_i\phi_j\\[1mm]
\|\nabla \phi_i\|^2= 1-\frac{1}{r^2}\phi_i^2\\[1mm]
\int_{\mathbb{S}^m_r}\phi^2_idM=\frac{r^2}{m+1}|\mathbb{S}^m_r|\\
\int_{\mathbb{S}^m}\|\nabla \phi_i\|^2dM
=\lambda_1(r)\int_{\mathbb{S}^2} \phi_i^2dM =\frac{m}{m+1}|\mathbb{S}^m_r|.
\end{array}\right.
\end{equation}
We also denote by $\int_{\mathbb{S}^m_r}\phi^2dM$ 
any of the integrals $\int_{\mathbb{S}^m_r}\phi_i^2dM$, $i=1,\ldots, m+1$.
\\[1mm]
Note that $\lambda_1(\mathbb{S}_r^m)=m\|H\|^2$,
$\|H\|=\frac{1}{r}$, and $|\mathbb{S}^m_r|=r^m|\mathbb{S}^m|$.
Any smooth function $f$ on $\mathbb{S}^m_r$ can be written as
a $L^2$-convergent sum $f=\sum_{l}\Psi_l$, where $\Psi_l=\sum_\sigma
A_{\sigma}\psi_{l,\sigma}$
is an $\lambda_l(r)$-eigenfunction and $A_{\sigma}$ are constants. 
This sum is in fact $H^1$-convergent
to $f$ ( see a proof of this in theorem 25.2 of \cite{Jost},  that
 formally holds for any compact Riemannian manifold as well).
If $l=1$ then $\phi^1,\ldots, \phi^{m+1}$ is up to a homothetic factor
an $L^2$-o.n. basis of
$E_{\lambda_1(r)}$.
If $\bar{f}$ is a homogeneous polynomial function of degree $l$ then,
for $X,Y\in T_xS^m_r$
$$\mathrm{Hess}\bar{f}(X,Y)=\mathrm{Hess}{f}(X,Y)-lf(x)\langle X,Y\rangle,$$
where $f$ is $\bar{f}$ restricted to $\mathbb{S}^m$.
If $r=1$,  the Ricci tensor of $\mathbb{S}^m$  is
$(m-1)\langle, \rangle$, and  using the Reilly's formula 
$$\mathrm{Ricci}(\nabla f, \nabla f)=(\Delta f)^2
+ \frac{1}{2}\Delta(\|\nabla f\|^2)
-\mathrm{div}(\Delta f \nabla f)-\|\mathrm{Hess}f\|^2$$
we obtain for $f\in E_{\lambda_l}$
$$ \begin{array}{l}
\lambda_l(\lambda_l-(m-1))=\frac{\int_{\mathbb{S}^{m}}
\|\mathrm{Hess}{f}\|^2dM}{\int_{\mathbb{S}^{m}}f^2dM}\\
\frac{(m-1)}{m}\lambda_l (\lambda_l- m)=
\frac{\int_{\mathbb{S}^{m}}\|\mathrm{Hess}{f}-\frac{\Delta f}{m}\langle,
\rangle\|^2dM}{\int_{\mathbb{S}^{m}}f^2dM}\neq 0~~~~~\mbox{if~}l\geq 2.
\end{array}$$
In particular $\mathrm{Hess}f$ is a multiple of the metric if and only if
$l=1$.

 The $\Omega$-stability condition for $\mathbb{S}^m_r$ is given
by  the inequality $I_{\Omega}(f\nu+ W,f\nu + W) \geq 0$
for any smooth section $W$ of $F$
and any smooth function  $f\in L_T^2(\mathbb{S}^m_r)$, where
$I_{\Omega}$ is given by  eq. (16) in \cite{Sa1},
\begin{eqnarray*}
I_{\Omega}(f\nu+ W,f\nu + W) &=& I(f,f)+m\|H\|\int_{\mathbb{S}_r^m}
f^2\bar{g}(C_{\Omega}(\nu),\nu)dM\\
&& +\int_{\mathbb{S}_r^m}
\left(\|\nabla^{\bot} W\|^2+m\|H\|\bar{g}(C_{\Omega}(W),W)\right) dM
\end{eqnarray*}
where 
$$I(f,f)=\int_{\mathbb{S}^m_r}\|\nabla f\|^2
-m\|H\|^2\int_{\mathbb{S}^m_r}f^2dM \geq 0$$ 
and $\bar{g}(C_{\Omega}(\nu),\nu)=\bar{\nabla}_{\nu}\Omega
(\nu, e_1, \ldots, e_m)$.
Then  if $C_{\Omega}\neq 0$, and
in particular if (8) does not hold we easily have an instability factor
for $\mathbb{S}_r^m$. For $f\in E_{\lambda_1(r)}$, we have $I(f,f)=0$ and 
$I_{\Omega}(f\nu,f\nu) \leq  r\lambda_1(r)b
\int_{\mathbb{S}_r^m} f^2dM$, 
where  
$$
b=\sup_{\mathbb{S}_r^m}\bar{\nabla}_{\nu}\Omega(\nu, e_1,\ldots, e_m).
$$
We also fix  a global parallel basis $W_{\alpha}$
 of $F=\mathbb{R}^{n-1}$, where $\alpha=m+2, \ldots m+n$.
Thus, by (2)  $$\int_{\mathbb{S}_r^m}
\bar{g}(C_{\Omega}(W_{\alpha}),W_{\alpha})dM=\int_{\mathbb{S}^m_r}
 \bar{\nabla}_{W_{\alpha}}\Omega(W_{\alpha},e_1, \ldots, e_m)dM.$$
Hence, we have the following conclusion:
\begin{proposition}
Suppose that $\mathbb{S}^m_r$ lies on a $\Omega$-calibrated vector
subspace $~\mathbb{R}^{m+1}$ of $~\mathbb{R}^{m+n}$ where $\Omega$ is a
semicalibration of rank $m+1$.
 If $b<0$ or if for some $\alpha\geq m+2$ we have
$\int_{\mathbb{S}^m_r}\bar{\nabla}_{W_{\alpha}}\Omega(W_{\alpha},
e_1, \ldots, e_m) dM<0$, then $\mathbb{S}^m_r$ is $\Omega$-unstable.
\end{proposition}

We now look for conditions for $\Omega$-stability to hold on spheres.
We define a $(m-1)$-form on $\mathbb{R}^{n+m}$ by
$$\hat{\xi}_{\alpha\beta}=\Omega(W_{\alpha},W_{\beta}, \ldots).$$
Then 
$$\xi(W_{\alpha},W_{\beta})=*\phi^*\hat{\xi}_{\alpha\beta} \in C^{\infty}(T^*M).$$
These forms are co-closed if $\Omega$ is parallel, but if
$\Omega$ is not parallel, the stability condition implies
co-closeness of $\xi(W_{\alpha},W_{\beta})$ as we will recall in
next theorem. 
\begin{theorem} [\cite{Sa1}] Let us suppose 
 $M=\mathbb{S}^m_r$ is a  $m$-sphere of radius $r$ of  a
$\Omega$-calibrated vector subspace $\mathbb{R}^{m+1}$ and that
$\Omega$ is a semicalibration such that (8) holds, that is
$\bar{\nabla}_{W}\Omega(W, e_1,\ldots, e_m)=0$ 
  for all $W\in NM$. \\[1mm]
Then $M$ is $\Omega$-stable if and only if the 1-forms 
$\xi(W_{\alpha}, W_{\beta})$ are co-exact, that is 
$$\xi(W_{\alpha}, W_{\beta})=\delta \omega_{\alpha\beta}$$
 for some 2-forms $\omega_{\alpha\beta}$ on $M$ and (10), or equivalently, (11)
 holds $\forall f_{\alpha}\in C^{\infty}(M)
,\alpha=m+2,\ldots, m+n $:\\[1mm]
\begin{eqnarray}
\sum_{\alpha<\beta}-2m\|H\|\int_{S^m_r} f_{\alpha}\xi(W_{\alpha},W_{\beta})(
\nabla f_{\beta}) &\leq& \sum_{\alpha}\int_M\|\nabla f_{\alpha}\|^2dM,\\[2mm]
\sum_{\alpha<\beta}-2 m\|H\|\int_{S^m_r} \omega_{\alpha\beta}(\nabla f_{\alpha},
\nabla f_{\beta})&\leq& \sum_{\alpha}\int_M\|\nabla f_{\alpha}\|^2dM.
\end{eqnarray} 
\end{theorem}
\noindent
The first inequality  (10) is a direct consequence of the
stability condition,  while (11) is proved in Proposition 4.5 of
\cite{Sa1} by using the  Hodge theory of spheres and that
$$\int_{\mathbb{S}^m_r}\xi(W_{\alpha},W_{\beta})(
f_{\alpha}\nabla f_{\beta}) dM=
\int_{\mathbb{S}^m_r}\langle \xi(W_{\alpha},W_{\beta}),
f_{\alpha}df_{\beta}\rangle dM=\int_{\mathbb{S}^m_r}
\omega_{\alpha\beta}(\nabla f_{\alpha},\nabla f_{\beta}) dM.$$
\begin{corollary} An $m$-sphere of an $\Omega$-calibrated Euclidean subspace
$\mathbb{R}^{m+1}$ of $\mathbb{R}^{m+n}$ for which $C_{\Omega}=0$
is $\Omega$-stable. This is 
the case when $n=2$ and $\Omega$ parallel. 
\end{corollary}
\noindent
The condition $C_{\Omega}=0$ is a very restrictive condition, and does not hold
for most calibrations coming from special holonomy, since (7) does not hold. 
But the operator $C_{\Omega}$ vanish
when  $\Omega$ is a semicalibration defined by a Riemannian fibration of 
$\mathbb{R}^{m+n}$ 
with some $(m+1)$-dimensional totally geodesic fibre where
$\mathbb{S}^m_r$ lies \cite{Sa1}.\\

Inequality (11) can be seen as the \em long integral Cauchy-Riemann 
$\Omega$-inequality \em for 
$(n-1)$-tuples of functions on $\mathbb{S}^m_r$.
If stability holds, then for each $\alpha<\beta$ fixed,
and  taking $f_{\gamma}=0$ for $\gamma\neq \alpha, \beta$,
 the  \em short integral Cauchy-Riemann $\Omega$-inequalities \em
hold for any pairs of functions $f,h$ on $\mathbb{S}^m_r$:
\begin{equation}
m\|H\|\left|\int_{S_r^{m}}\omega_{\alpha\beta}(\nabla f, \nabla h)dM\right|
\leq  \|\nabla f\|_{L^2}\|\nabla h\|_{L^2} ~~~\forall \alpha<\beta.
\end{equation}
On the other hand, inequality (10) gives us a tool to determine
if a sphere is $\Omega$-stable applying the Rayleigh characterization
of the spectrum of $\mathbb{S}^m_r$. 
Let $\Theta(r)=\sup_{\alpha<\beta}\Theta_{\alpha\beta}(r)$  where
$\Theta_{\alpha\beta}(r)=\sup_{\mathbb{S}^m_r}
\|\xi(W_{\alpha},W_{\beta})\|\leq 1$. 
For $f_{\alpha}\in E_{\lambda_{l_{\alpha}}(r)}^+$, by Schwartz inequality,
\begin{equation}
\left|\int_{\mathbb{S}^m_r}
f_{\alpha}\xi(W_{\alpha},W_{\beta})(\nabla f_{\beta})\, dM\right|
\leq \Theta_{\alpha\beta}(r) \frac{1}{\sqrt{\lambda_{l_{\alpha}}(r)}}
\|\nabla f_{\alpha}\|_{L^2}\|\nabla f_{\beta}\|_{L^2}.
\end{equation}
Using the inequality 
$\sum_{m+2\leq \alpha<\beta\leq m+n}2a_{\alpha}a_{\beta}\leq (n-2)
(a_{m+2}^2+ \ldots + a_{m+n}^2)$, for any real numbers
$a_{\alpha}$,
and  $\lambda_l(r)=\frac{\lambda_l}{r^2}$  we get  next Proposition:
\begin{proposition} Assuming $f_{\alpha}\in E_{\lambda_{l}(r)}^+$
for all $\alpha$ then
$$-\sum_{\alpha<\beta}{2m}\|H\|\int_{S^{m}_r} f_{\alpha}\xi(W_{\alpha},W_{\beta})(\nabla f_{\beta})\, dM~\leq ~\, 
\frac{m(n-2)}{\sqrt{\lambda_l}}\Theta(r)\sum_{\alpha}
\int_{\mathbb{S}^m_r}\|\nabla f_{\alpha}\|^2dM.$$
Consequently, $\mathbb{S}_r^m$ is $\Omega$-stable in $\mathbb{R}^{m+n}$
if  $n=2$, or if $n\geq 3$ and $\Theta(r)\leq \frac{1}{\sqrt{m}(n-2)}$.
\\[1mm]
\end{proposition}
\begin{corollary} Supposing that $\xi(W_{\alpha},W_{\beta})=\delta
\omega_{\alpha\beta}$,  
and $f_{\alpha}\in E_{\lambda_l(r)}^{+}$ for all $\alpha$, then 
$$\sum_{\alpha<\beta}-{2m}\|H\|\int_{S^{m}_r} \omega_{\alpha\beta}
(\nabla f_{\alpha},\nabla f_{\beta})dM~\leq ~\, 
\frac{m(n-2)}{\sqrt{\lambda_l}}\Theta(r)\sum_{\alpha}
\int_{\mathbb{S}^m_r}\|\nabla f_{\alpha}\|^2dM.$$
Hence, if $l$ is sufficiently large such that  $l(l+m-1)\geq m^2(n-2)^2$,  
the long Cauchy Riemann inequality (11) holds
for functions $f_{\alpha}\in E_{\lambda_l(r)}^+$.
\end{corollary}
\noindent
This  estimate is in general not sharp, because it ignores the
signs that $\omega_{\alpha\beta}$ can take.
We also remark that for $x\in \mathbb{S}^{m}$, since $T_x \mathbb{S}^{m}
=T_{rx} \mathbb{S}^{m}_r$, if $\xi(x)=\xi(rx)$ then
 the same holds for for $\omega_{\alpha\beta}$ and
$\Theta_{\alpha\beta}(r)=\Theta_{\alpha\beta}$.
In this case  $\mathbb{S}^{m}$ is $\Omega$-stable if and only if
 $\mathbb{S}^{m}_r$ is so.

If $\xi(W_{\alpha},W_{\beta})=\delta \omega_{\alpha\beta}$
then
$$\int_{\mathbb{S}^m_r}\omega_{\alpha\beta}(\nabla f_{\alpha},
\nabla f_{\beta})dM =\frac{1}{2}\int_{\mathbb{S}^m_r}\La{(}
f_{\alpha}\xi(W_{\alpha},W_{\beta})(\nabla f_{\beta})-
f_{\beta}\xi(W_{\alpha},W_{\beta})(\nabla f_{\alpha})\La{)}dM.$$
Applying  inequality (13) to this expression 
we immediately deduce that:
\begin{proposition} 
If we fix  
$\alpha<\beta$ and $\xi(W_{\alpha},W_{\beta})=\delta \omega_{\alpha\beta}$, 
then for  any functions 
$f\in E_{\lambda_l(r)}^+$ and $h\in E_{\lambda_s(r)}^+$
we  have 
$$2{m}\|H\|\left|\int_{\mathbb{S}^m_r} \omega_{\alpha\beta}
(\nabla f,\nabla h)dM
\right|\leq {m}\Theta_{\alpha\beta}(r)
\left(\frac{1}{\sqrt{\lambda_{l}(r)}}
+\frac{1}{\sqrt{\lambda_{s}(r)}}\right)
\|\nabla f\|_{L^2}\|\nabla h\|_{L^2}.$$
\end{proposition}
\section{The 2-sphere of $\mathbb{R}^7$}
In this section we consider the unit 2-sphere $\mathbb{S}^2$ of an associative
Euclidean 3-dimensional subspace $\mathbb{R}^3$ of $\mathbb{R}^7$, that
is, we may assume $\mathbb{R}^3=span\{\epsilon_1, \epsilon_2, \epsilon_3\}$
is $\Omega$-calibrated by the associative calibration $\Omega$
defined in the introduction.
As we have  pointed out in Remark 4.4 \cite{Sa1}, taking
$W_{\alpha}=\epsilon_{\alpha}$ for $\alpha=4,5,6,7$,  then
$$\hat{\xi}_{45}=\hat{\xi}_{67}=\epsilon_*^1=dx^1~~~~~
\hat{\xi}_{46}=-\hat{\xi}_{57}=\epsilon_*^2=dx^2~~~~~
\hat{\xi}_{47}=\hat{\xi}_{56}=-\epsilon_*^3=-dx^3.$$
Since $\delta\omega_{\alpha\beta}=
\xi(W_{\alpha},W_{\beta})=*\phi^*\hat{\xi}_{\alpha\beta}$,  and 
$\omega_{\alpha\beta}=\rho_{\alpha\beta}Vol_{\mathbb{S}^2}$
we conclude that
$$\rho_{45}=\rho_{67}=-\phi_1~~~~~
\rho_{46}=-\rho_{57}=-\phi_2~~~~~
\rho_{47}=\rho_{56}=\phi_3$$
(there is a misprint in \cite{Sa1}, a wrong sign for $\hat{\xi}_{56}$).
Note that $Vol_{\mathbb{S}^2}(X,Y)=\langle JX,Y\rangle$.
Then (3) holds iff (11) holds, and (4) holds iff (12) holds. We first prove  
 Theorem 1.1.
We will need some lemmas:
\begin{lemma} On a K\"{a}hler manifold $(M,J,g)$ of real dimension $2k$, for any
functions $ f,h\in C^{\infty}(M)$ we have
$$ g( J\nabla f,\nabla h) =\mathrm{div}\la{(}
h\, J\nabla f\la{)}=
\sm{\frac{1}{2}}\,\mathrm{div}\la{(}
h\, J\nabla f-f\, J\nabla h\la{)}.$$
Furthermore, if $M$ is a closed manifold, then  the operator
$$\eta(\phi, f,h):=\int_M\phi g( J\nabla f, \nabla h) dM$$
is skew-symmetric in the three variables $\phi, f, h \in C^{\infty}(M)$
\end{lemma}
\noindent
\em Proof. \em Let $e_i, i=1, \ldots, 2k$ be a local o.n. frame 
of $TM$ with $e_{k+i}=Je_i$, $i=1, \ldots, k$. Then
\begin{eqnarray*}
\mathrm{div}\la{(}
h\, J\nabla f\la{)}&=& \sum_{1\leq i\leq 2k}g(\nabla_{e_i}(hJ\nabla f), e_i)
=\sum_{1\leq i\leq 2k}dh(e_i)g(J\nabla f, e_i)-h\, Hess\, f(Je_i, e_i)\\
&=& g( J\nabla f, \nabla h)-\!\!\!\!\sum_{1\leq i\leq k}\!\!\!\!\La{(}
h\, Hess\, f(Je_i, e_i)+h\,   Hess\, f(JJe_i,J e_i)\La{)}
=g( J\nabla f, \nabla h)
\end{eqnarray*}
where  we used the fact that $\mathrm{Hess}f$ is symmetric.
The second equality follows immediately.
Using the equality
$$ \mathrm{div}(\phi h J\nabla f)=\phi\, \mathrm{div}(hJ\nabla f)+
g(\nabla \phi, h J\nabla f)
=\phi\, \mathrm{div}(hJ\nabla f)+h\, \mathrm{div}(\phi J\nabla f)$$
and applying Stokes we see that $\eta$ is skew-symmetric.
\qed\\[3mm]
Next lemma is a Weitzenb\"{o}ck type formula:
\begin{lemma} On a K\"{a}hler manifold $(M,J,g)$ of real dimension $2k$
we have for any functions $f,h$ and o.n. frame $e_i$, $i=1, \ldots, 2k$
\begin{eqnarray*}
\Delta \left(g(J\nabla f, \nabla h)\right) &=& g(J\nabla \Delta f, \nabla h)
+g(J\nabla f, \nabla \Delta h)-\mathrm{Ricci}(\nabla f, J\nabla h)
+\mathrm{Ricci}(\nabla h, J\nabla f)\\
&&-\sum_{1\leq i,j\leq 2k}
2\, \mathrm{Hess}\, f(e_i,Je_j)\mathrm{Hess}\ h(e_i,e_j).
\end{eqnarray*}
\end{lemma}
\noindent
\em Proof. \em We may assume at a point $p$, $\nabla e_i(p)=0$.
Differentiating 
$d(g(J\nabla f,\nabla h))(e_i)=g(J\nabla_{e_i}\nabla f, \nabla h)
+g(J\nabla f, \nabla_{e_i}\nabla h)$ with respect to $e_i$ we have
at the point $p$
\begin{eqnarray*}
\lefteqn{\Delta \left(g(J\nabla f, \nabla h)\right)=}\\
&&=\sum_{i} g(J\nabla_{e_i}\nabla_{e_i}\nabla f, \nabla h)
+ 2g(J\nabla_{e_i}\nabla f, \nabla_{e_i}\nabla  h)+
g(J\nabla f, \nabla_{e_i}\nabla_{e_i}\nabla h)\\
&=& \sum_{ij}\nabla^2_{e_i, e_i}d f(e_j) d h(Je_j)
-2\mathrm{Hess}\, f(e_i, Je_j)\mathrm{Hess}\, h (e_i,e_j)
-df(Je_j)\nabla^2_{e_i,e_i}dh(e_j)
\end{eqnarray*}
where $\nabla^2_{X,Y}df=\nabla_{X}\nabla_{Y}df-\nabla_{\nabla_XY}df$,
for any vector fields $X,Y,Z$.
Here we use the curvature sign
$R(X,Y)=-\nabla_{X}\nabla_Y+\nabla_Y\nabla_X + \nabla_{[X,Y]}$.
Then we have (see e.g.\  \cite{Sa2} p.1234)
$$\nabla^2_{X,Y}df(Z)=\nabla^2_{X,Z}df(Y)=\nabla^2_{Y,X}df(Z)+ df(R^M(X,Y)Z).$$
Thus
$$
\sum_i\nabla^2_{e_i, e_i}d f(e_j) 
=\sum_i\nabla^2_{e_j, e_i}d f(e_i)
+df(R^M(e_i,e_j)e_i)
=\nabla_{e_j}(\Delta f)+ df(Ricci^M(e_j)).$$
Replacing this equality in the above equation, and a similar one w.r.t. $h$
we obtain the formula  of the lemma.\qed
\begin{lemma} $(1)$~If $f\in E_{\lambda_1}$ and $h\in E_{\lambda_r}$ then
$\langle J\nabla f,\nabla h\rangle \in E_{\lambda_r}$.\\[2mm]
$(2)$ $\langle J\nabla \phi_i,\nabla \phi_j\rangle =\phi_k$~
for $(i,j,k)$ is a positive permutation of $(1,2,3)$.\\[2mm]
$(3)$ If $f\in E_{\lambda_l}$ and $h\in E_{\lambda_r}$,
then $$\Delta \langle J\nabla f, \nabla h\rangle =-(\lambda_l+\lambda_r-2)
\langle J\nabla f, \nabla h\rangle-\sum_{ij} 2\mathrm{Hess}\, f(e_i,Je_j)
\mathrm{Hess}\, h(e_i,e_j).$$
\end{lemma}
\noindent
\em Proof. \em (1) Since $\phi_1, \phi_2, \phi_3$ is a basis of 
$\lambda_1$-eigenfunctions, then $\mathrm{Hess}\,f=-f\cdot g$. 
Moreover $\mathrm{Ricci}^M=g$.
Taking $e_i$ a o.n. basis that diagonalizes $\mathrm{Hess} \, h$ we conclude 
from Lemma 3.2 that
$\langle J\nabla f,\nabla h\rangle \in E_{\lambda_r}$.
(2) If we consider spherical coordinates $
\phi:[0,\pi]\times[0,2\pi]\to \mathbb{S}^2$,~
$\phi(\varphi,\theta)=(\sin\varphi\cos\theta,\sin\varphi\sin\theta,
\cos\varphi)$,
then $X=\frac{\partial}{\partial \varphi}$, $Y=\frac{1}{\sin\varphi}
\frac{\partial}{\partial \theta}$ defines a d.o.n. frame of $T\mathbb{S}^2$
and so $JX=Y$. Furthermore, 
$$\nabla \phi_1=\cos\varphi\cos\theta \, X-\sin\theta\, Y,~~~~
\nabla \phi_2=\cos\varphi\sin\theta \, X+ \cos\theta\, Y,~~
\nabla \phi_3=-\sin\varphi\, X$$
Then we see that (2) holds. (3) is an immediate consequence of  Lemma 3.2
and generalizes (1).
\qed
\begin{proposition}
$(1)$~ If $f\in E_{\lambda_l}$ and $h\in E_{\lambda_r}$ with $l\neq r$
then $\int_{\mathbb{S}^2}\phi_i\langle J\nabla f, \nabla h\rangle dM=0$.
\\[2mm]
$(2)$~ $2\int_{\mathbb{S}^2}\phi_k\langle J\nabla \phi_i,\nabla \phi_j\rangle
=\epsilon_{kij}\|\nabla\phi_i\|_{L^2} \|\nabla \phi_j\|_{L^2}$,
where $\epsilon_{ikj}=+1, -1$, according to the signature
of $(i,j,k)$ as a permutation of $(1,2,3)$, or zero if repeated indexes
appear. 
\end{proposition}
\noindent
\em Proof. \em  Using Lemma 3.1 
$$\int_{\mathbb{S}^2}\phi_i\langle J\nabla f, \nabla h\rangle =
-\int_{\mathbb{S}^2}h \langle J\nabla \phi_i, \nabla f\rangle dM.$$
By Lemma 3.3 (1) $\langle J\nabla \phi_i, \nabla f\rangle\in E_{\lambda_l}$,
and so this is $L^2$-orthogonal to $h\in E_{\lambda_r}$.\\
(2) is an immediate consequence of Lemma 3.3 (2) and (9). \qed\\[4mm]
\noindent
{\bf Proof of Theorem 1.1.}\\[1mm]
 We may assume $i=1$.
Inequality (4) is equivalent to prove that 
$$4\int_{\mathbb{S}^2}\phi_1\langle J\nabla f, \nabla h
\rangle dM\leq \|\nabla f\|_{L^2}^2+ \|\nabla h\|^2_{L^2},$$
holds  $\forall f, h\in
C^{\infty}(\mathbb{S}^2)$.
We write $f=f_0+f_1+ f'$ and $h=h_0+h_1+ h'$ where $f_0,h_0 \in E_{\lambda_0}$
 are constants,
$f_1, h_1\in E_{\lambda_1}$, and $f', h'\in E_{\lambda_2}^+$.
Then applying Proposition 3.1 (1)
$$4\int_{\mathbb{S}^2}\phi_1\langle J\nabla f, \nabla h\rangle dM
=4\int_{\mathbb{S}^2}\phi_1\langle J\nabla f_1, \nabla h_1\rangle dM
+4\int_{\mathbb{S}^2}\phi_1\langle J\nabla f', \nabla h'\rangle dM$$
From Proposition 2.3 we have
$$4\int_{\mathbb{S}^2}\phi_1\langle J\nabla f', \nabla h'\rangle dM
\leq \frac{4}{\sqrt{6}}\|\nabla f'\|_{L^2}\|\nabla h'\|_{L^2}
\leq \frac{2}{\sqrt{6}}(\|\nabla f'\|_{L^2}+\|\nabla h'\|_{L^2}).$$
Since $f_1=\sum_i \mu_i\phi_i$, $h_1=\sum_j \sigma_j\phi_j$
where $\mu_i,\sigma_j$ are constants, 
applying Proposition 3.1(2) and (9) we have 
\begin{eqnarray*}
4\int_{\mathbb{S}^2}\phi_1\langle J\nabla f_1, \nabla h_1\rangle dM
&=&4(\mu_2\sigma_3-\mu_3\sigma_2)\int_{\mathbb{S}^2}\phi_1^2dM\\
&\leq& (\sum_i \mu_i^2+\sigma_i^2)\int_{\mathbb{S}^2}\|\nabla \phi_1\|^2dM
=\|\nabla f_1\|^2_{L^2}+\|\nabla h_1\|^2_{L^2}.
\end{eqnarray*}
As $\frac{2}{\sqrt{6}}<1$ we conclude  that (4) holds, and equality
is achieved if and only if $f'=h'=0$ and $\mu_1=\sigma_1=0$,
$\mu_2=\sigma_3$, $\mu_3=-\sigma_2$. \qed\\[3mm]

Next we prove Theorem 1.2. in several steps, as a consequence of
Proposition 3.3, and of Lemmas 3.7, 3.8, and 3.10. Indeed, 
as a consequence of Proposition 3.1(1), 
we only need to droop our attention on $4$-tuples 
$(f_4,f_5,f_6,f_7)$ that have at least two components in the same eigenspace.
The case we have two pairs of functions in two different eigenspaces
lies in the case of Theorem 1.1 as we can easily verify. Thus we have
\begin{proposition} If two  elements of $\{f_4,f_5,f_6,f_7\}$
are in $E_{\lambda_r}$ and the other two in $E_{\lambda_l}$
with $r\neq l$, then the long Cauchy-Riemann inequality (3) holds.
\end{proposition}
\noindent
If only three of the functions $f_{\alpha}$
 are in the same eigenspace, the terms involving
the forth function vanish, and so we may assume the later to be zero, 
that is we are in the case that all functions are in the same eigenspace.
This is the case we are now considering.\\

We denote by $\epsilon_i$, $i=1,2,3$  the canonical
 basis of $\mathbb{R}^3$,
and so  a multi-index of nonnegative integers 
is of the form $a=  (a_1,a_2,a_3)=\sum_{i}a_i\epsilon_i$.
Next we recall the well known formula ( see for instance \cite{Bre} appendix)
\begin{lemma} If $P:\mathbb{R}^3\to \mathbb{R}$ is an homogeneous polynomial
function
of degree $l$ then
$$\int_{\mathbb{S}^2}P(x)dM=
\frac{1}{\lambda_l}\int_{\mathbb{S}^2}\Delta^0 P(x)dM.$$
In particular for $|a|=a_1+a_2+a_3=l$
$$\int_{\mathbb{S}^2}\phi^{a}dM=\sum_{1\leq i\leq 3}\frac{a_i
(a_i-1)}{l(l+1)}
\int_{\mathbb{S}^2}\phi^{a-2\epsilon_i}dM,$$
where the terms  with $a_i<2$ are considered to vanish.
\end{lemma}
\noindent
Let us denote by $(O)$ and $(E)$ meaning odd and even respectively.
We also say that $a$ is $(O,E,E)$  meaning that $a_1$ is odd and
$a_2$ and $a_3$ are even, and so on. 
Since $\int_{\mathbb{S}^2}\phi_idM=0$, by induction we conclude
from the previous lemma that 
$$\int_{\mathbb{S}^2}\phi^{a}dM=0~~~\mbox{iff}~~\exists i: a_i\mbox{~is
(O)}.$$
Now using (9) and Lemma 3.3(2) 
we obtain the  following two lemmas, respectively:
\begin{lemma} If $|a|=|b|=l$ then
\begin{eqnarray*}
\int_{\mathbb{S}^2}\langle \nabla \phi^{a}, \nabla \phi^{b}\rangle dM
&=&- l^2\int_{\mathbb{S}^{2}}
\phi^{a+b}dM
+ \sum_i a_i b_i\int_{\mathbb{S}^{2}}\phi^{a+b-2\epsilon_i}dM\\
&=&\sum_i\frac{l((a_i+b_i)-(a_i-b_i)^2)+2a_ib_i}{2(2l+1)}
\int_{\mathbb{S}^{2}}\phi^{a+b-2\epsilon_i}dM. \nonumber
\end{eqnarray*}
If this does not vanish then $a+b$ is $(E,E,E)$.
\end{lemma}
\begin{lemma} 
$$\begin{array}{llcl}
(1)~~~~~~
&\int_{\mathbb{S}^2}\phi_1\langle J\nabla\phi^{a},\nabla \phi^{b}\rangle
dM &=&(a_1b_2-a_2b_1)\int_{\mathbb{S}^2}
\phi^{a+b-\epsilon_2+\epsilon_3}dM\\[1mm]
~~& &&+(-a_1b_3+a_3b_1)\int_{\mathbb{S}^2}
\phi^{a+b+\epsilon_2-\epsilon_3}dM\\[1mm]
~~& &&+(a_2b_3-a_3b_2)\int_{\mathbb{S}^2}
\phi^{a+b +2\epsilon_1-\epsilon_2-\epsilon_3}dM.
\end{array}$$
If this does not vanish then $a+b$ is $(E,O,O)$.\\[1mm]
$$\begin{array}{llcl}
(2)~~~~~~~
&\int_{\mathbb{S}^2}\phi_2\langle J\nabla\phi^{a},\nabla \phi^{b}\rangle
dM &=&(a_2b_3-a_3b_2)\int_{\mathbb{S}^2}
\phi^{a+b+\epsilon_1-\epsilon_3}dM\\[1mm]
~~& &&+(-a_2b_1+a_1b_2)\int_{\mathbb{S}^2}
\phi^{a+b-\epsilon_1+\epsilon_3}dM\\[1mm]
~~& &&+(a_3b_1-a_1b_3)\int_{\mathbb{S}^2}
\phi^{a+b -\epsilon_1+ 2\epsilon_2-\epsilon_3}dM.
\end{array}$$
If this  does not vanish  then $a+b$ is $(O,E,O)$.\\[1mm]
$$\begin{array}{llcl}
(3)~~~~~~~
&\int_{\mathbb{S}^2}\phi_3\langle J\nabla\phi^{a},\nabla \phi^{b}\rangle
dM &=&(a_3b_1-a_1b_3)\int_{\mathbb{S}^2}
\phi^{a+b-\epsilon_1+\epsilon_2}dM\\[1mm]
~~& &&+(-a_3b_2+a_2b_3)\int_{\mathbb{S}^2}
\phi^{a+b+\epsilon_1-\epsilon_2}dM\\[1mm]
~~& &&+(a_1b_2-a_2b_1)\int_{\mathbb{S}^2}
\phi^{a+b -\epsilon_1-\epsilon_2+ 2\epsilon_3}dM.\\[2mm]
\end{array}$$
If this does not vanish then $a+b$ is $(O,O,E)$.\\
\end{lemma}
\noindent
Note that since $|a|=|b|= l$ then
$$
\begin{array}{l}
-a_1b_3+a_3b_1 = a_1b_2-a_2b_1 + l(b_1-a_1)\\
a_2b_3-a_3b_2 =a_1b_2-a_2b_1 + l(a_2-b_2).
\end{array}
$$
\begin{proposition} If $f_4, f_5,f_6,f_7$  are all elements of $E_{\lambda_l}$
where $l\geq 6$, then the long Cauchy-Riemann inequality (3) holds.
\end{proposition}
\noindent
\em Proof. \em  We have $\Theta(1)=1$ and $m=2$, $n=5$. Therefore,  
$m^2(n-2)^2=36\leq \lambda_6= 6\times 7$.
From Corollary 2.2.  we conclude that (3) holds for $l\geq 6$. \qed\\[4mm]

We now have to consider the case
$l\leq 5$.
We write $f_{\alpha}=\sum_{|a|=l}A^{\alpha}_{a}\phi^{a}$,
where $A^{\alpha}_{a}$ are constants.
This summation is not $L^2$-orthogonal, in general,  since  eigenfunctions
$\psi_{l,\sigma}$ are  usually a sum of linearly  independent monomials.
We have using Lemma 3.6
\begin{eqnarray}
\lefteqn{\begin{array}{l}
4\int_{\mathbb{S}^2}\phi_1 (\langle J\nabla f_4,\nabla f_5\rangle
+\langle J\nabla f_6,\nabla f_7\rangle )dM \\[1mm]
+ 4\int_{\mathbb{S}^2}\phi_2 (\langle J\nabla f_4,\nabla f_6\rangle
-\langle J\nabla f_5,\nabla f_7\rangle )dM  \\[1mm]
- 4\int_{\mathbb{S}^2}\phi_3 (\langle J\nabla f_4,\nabla f_7\rangle
+\langle J\nabla f_5,\nabla f_6\rangle )dM =\end{array}}\\
&=& \!\!\!\!\!\!\sum_{a+b=(E,O,O)}\!\!\!\!\!\!\!\!\!\! 
4(A^{4}_{a}A^{5}_{b}+A^{6}_{a}A^{7}_{b})
\int_{\mathbb{S}^2}\phi_1\langle J\nabla \phi^a,\nabla\phi^b\rangle dM\\
&& +\!\!\!\!\!\!\sum_{a+b=(O,E,O)} \!\!\!\!\!\!\!\!\!
4(A^{4}_{a}A^{6}_{b}-A^{5}_{a}A^{7}_{b})
\int_{\mathbb{S}^2}\phi_2\langle J\nabla \phi^a,\nabla\phi^b\rangle dM\\ 
&& -\!\!\!\!\!\!\sum_{a+b=(O,O,E)} \!\!\!\!\!\!\!\!\!
4(A^{5}_{a}A^{6}_{b}+A^{4}_{a}A^{7}_{b})
\int_{\mathbb{S}^2}\phi_3\langle J\nabla \phi^a,\nabla\phi^b\rangle dM.
\end{eqnarray} 
We will divide  the proof into several lemmas.
\begin{lemma} If $l=1$ and set $A^{\alpha}_{ i}
=A^{\alpha}_{\epsilon_i}$, that is, 
$f_{\alpha}= A^{\alpha}_{1}\phi_1+ A^{\alpha}_{2}\phi_2 
+A^{\alpha}_{3}\phi_3$,
then (4) holds. Furthermore, given 
$f_4$ and $f_5$ with  arbitrary coefficients $A^4_{i}$ and $A^5_{i}$
respectively,  then equality holds
for a 4-tuple $(f_4,f_5,f_6, f_4)$
 iff
$$\begin{array}{l}
f_6=(A^5_{2}+A^4_{3})\phi_1+A^6_{2}\phi_2+A^6_{3}\phi_3\\
f_{7}=(A^4_{2}-A^5_{3})\phi_1-(A^4_{1}+A^6_{3})\phi_2 
+(A^5_{1}+A^6_{2})\phi_3.
\end{array}$$
\end{lemma}
\noindent
\em Proof. \em Using Lemmas 3.6 and 3.3(2)
\begin{eqnarray*}
(15)+(16)+(17) &=& \La{(}4(A^4_{2}A^5_{3}-A^4_{3}A^5_{2}+A^6_{2}A^7_{3}-A^6_{3}A^7_{2})\\
&&\quad -4(A^4_{1}A^6_{3}-A^4_{3}A^6_{1}-A^5_{1}A^7_{3}+A^5_{3}A^7_{1})\\
&&\quad -4(A^4_{1}A^7_{2}-A^4_{2}A^7_{1}+A^5_{1}A^6_{2}-A^5_{2}A^6_{1})\La{)}
\int_{\mathbb{S}^2}\phi^2dM.
\end{eqnarray*}
Since $(-a+b+c)^2\geq 0$ and $(a+b+c)^2\geq 0$
 for any real numbers $a,b,c$,  we have
\begin{equation}
\left\{\begin{array}{ll}
2ac+2ab-2bc \leq  a^2+b^2+c^2~~~&\mbox{with equality iff }a=b+c\\
-2ac-2ab-2bc\leq a^2+b^2+c^2~~~&\mbox{with equality iff }a=-b-c \end{array}
\right.
\end{equation}
Applying these inequalities,  we have
\begin{eqnarray*}
\lefteqn{(15)+(16)+(17)=
 \La{(} 4(A^4_{2}A^5_{3}+A^4_{2}A^7_{1}-A^5_{3}A^7_{1})
+4(-A^4_{3}A^5_{2}+A^4_{3}A^6_{1}-A^5_{2}A^6_{1})}\\
&&~~ +4(A^6_{2}A^7_{3}+A^5_{1}A^7_{3}-A^5_{1}A^6_{2})
+4(-A^6_{3}A^7_{2}-A^4_{1}A^6_{3}-A^4_{2}A^7_{2}) \La{)}
\int_{\mathbb{S}^2}\phi^2dM\\
&\leq & \La{(} (A^4_{2})^2+(A^5_{3})^2+(A^7_{1})^2
+ (A^4_{3})^2+(A^5_{2})^2+ \\
& & ~~+(A^6_{1})^2 +(A^6_{2})^2+(A^7_{3})^2+(A^5_{1})^2
+(A^6_{3})^2+(A^7_{2})^2+ (A^4_{1})^2\La{)}
\int_{\mathbb{S}^2}\phi^2dM,
\end{eqnarray*}
with equality iff 
\begin{equation}
\left\{\begin{array}{lcl}
A^4_{2}=A^5_{3}+A^7_{1}&\quad&
A^6_{1}=A^4_{3}+A^5_{2}\\
A^7_{3}= A^5_{1}+A^6_{2}&\quad&
A^7_{2}=-A^4_{1}-A^6_{3}.\end{array}\right.
\end{equation}
On the other hand, using the last equality of (9) with $\lambda_1=2$ we see that
$$
\sum_{\alpha}\|\nabla f_{\alpha}\|^2_{L^2}
=(\sum_{\alpha i} 2(A^{\alpha}_{ i})^2)
\int_{\mathbb{S}^2}\phi^2dM.
$$
Thus,  $(15)+(16)+(17)\leq \sum_{\alpha}\|\nabla f_{\alpha}\|^2_{L^2}$ 
with equality iff (19) holds, and the lemma is proved. \qed\\[5mm]
Next we consider the case $l\geq 2$.
{\small
\begin{eqnarray}
\lefteqn{(15)=\sum_{a+b=(E,O,O)}
4(A^4_{a}A^5_{b}+A^6_{a}A^7_{b})\cdot}\\
&&\cdot\left(\ft{(a_1b_2-a_2b_1)\int_{\mathbb{S}^2}
\phi^{a+b-\epsilon_2+\epsilon_3}
+(-a_1b_3+a_3b_1)\int_{\mathbb{S}^2}
\phi^{a+b+\epsilon_2-\epsilon_3}
+(a_2b_3-a_3b_2)\int_{\mathbb{S}^2}
\phi^{a+b +2\epsilon_1-\epsilon_2-\epsilon_3}}\right)\nonumber
\end{eqnarray}
\begin{eqnarray}
\lefteqn{(16)=\sum_{a+b=(O,E,O)}
4(A^4_{a}A^6_{b}-A^5_{a}A^7_{b})\cdot}\\
&&\cdot  \left(\ft{(a_2b_3-a_3b_2)\int_{\mathbb{S}^2}
\phi^{a+b+\epsilon_1-\epsilon_3}
+(-a_2b_1+a_1b_2)\int_{\mathbb{S}^2}
\phi^{a+b-\epsilon_1+\epsilon_3}
+(a_3b_1-a_1b_3)\int_{\mathbb{S}^2}
\phi^{a+b -\epsilon_1+ 2\epsilon_2-\epsilon_3}}\right)\nonumber
\end{eqnarray}
\begin{eqnarray}
\lefteqn{(17)=\sum_{a+b=(O,O,E)}
4(-A^5_{a}A^6_{b}-A^4_{a}A^7_{b})\cdot }\\
&&\cdot \left(\ft{ (a_3b_1-a_1b_3)\int_{\mathbb{S}^2}
\phi^{a+b-\epsilon_1+\epsilon_2}
+(-a_3b_2+a_2b_3)\int_{\mathbb{S}^2}
\phi^{a+b+\epsilon_1-\epsilon_2}
+(a_1b_2-a_2b_1)\int_{\mathbb{S}^2}
\phi^{a+b -\epsilon_1-\epsilon_2+ 2\epsilon_3}}\right)\nonumber
\end{eqnarray}}
Note that the above terms are such that $a\neq b$ ( otherwise $a+b=(E,E,E)$),
and  are skew-symmetric on  $(a,b)$. 
We also have
\begin{equation}
\sum_{\alpha}\|\nabla f_{\alpha}\|^2 = 
\sum_{ab}\left[ (\sum_{\alpha}A_{\alpha a}A_{\alpha b})
\left(\sum_i\frac{l((a_i+b_i)-(a_i-b_i)^2)+2a_ib_i}{2(2l+1)}
\int_{\mathbb{S}^{2}}\phi^{a+b-2\epsilon_i}\right)\right]
\end{equation}
\begin{lemma} If $l=2$, then (4) holds with equality if and only if
$f_{\alpha}=0$ ~~$\forall \alpha$.
\end{lemma}
\noindent
\em Proof. \em
So we have $a,b$ running over 
$$\{(2,0,0),(0,2,0),(0,0,2),(1,1,0),(1,0,1),
(0,1,1)\}.$$
In $(20)=(15)$ we have to consider the following  terms: 
The terms with $a+b=(2,1,1)$ that are given by 
$a=(2,0,0)$ with $b=(0,1,1)$,
$a=(1,1,0)$ with $b=(1,0,1)$, and vice versa. 
The terms with $a+b=(0,3,1)$ that are given by
$a=(0,2,0)$ with $b=(0,1,1)$ and vice versa.
The terms with $a+b= (0,1,3)$ that are given by
$a=(0,0,2)$ with $b=(0,1,1)$ and vice versa.
Thus,
{\small 
\begin{eqnarray*}
\lefteqn{(15)=}\\
&&\begin{array}{l}
=  4 (A^4_{(200)}A^5_{(011)}-A^4_{(011)}A^5_{(200)}+
A^6_{(200)}A^7_{(011)}-A^6_{(011)}A^7_{(200)})
\cdot \La{(} 2\int_{\mathbb{S}^2}\phi_1^2\phi_3^2
-2\int_{\mathbb{S}^2}\phi_1^2\phi_2^2\La{)}\\
+  4 (A^4_{(110)}A^5_{(101)}-A^4_{4(101)}A^5_{(110)}+
A^6_{(110)}A^7_{(101)}-A^6_{(101)}A^7_{(110)})\cdot
 \La{(} -\int_{\mathbb{S}^2}\phi_1^2\phi_3^2
-\int_{\mathbb{S}^2}\phi_1^2\phi_2^2 
+ \int_{\mathbb{S}^2}\phi_1^4 \La{)}\\
+  4 (A^4_{(020)}A^5_{(011)}-A^4_{(011)}A^5_{(020)}+
A^6_{(020)}A^7_{(011)}-A^6_{(011)}A^7_{(020)})\cdot
 \La{(} 2\int_{\mathbb{S}^2}\phi_1^2\phi_2^2\La{)}\\
+  4 (A^4_{(002)}A^5_{(011)}-A^4_{(011)}A^5_{(002)}+
A^6_{(002)}A^7_{(011)}-A^6_{(011)}A^7_{(002)})
\cdot \La{(} -2\int_{\mathbb{S}^2}\phi_1^2\phi_3^2 \La{)}.\\
\end{array}
\end{eqnarray*} }
Using lemma 3.4, we have
$\int_{\mathbb{S}^2}\phi_i^2\phi_j^2dM=\frac{1}{5}(1+2\delta_{ij})
\int_{\mathbb{S}^2}\phi^2dM$.
Therefore,
{\small
\begin{eqnarray*}
\begin{array}{ll}
(15)=& \left[\frac{4}{5}\la{(}A^4_{(110)}A^5_{(101)}-A^4_{(101)}A^5_{(110)}+
A^6_{(110)}A^7_{(101)}-A^6_{(101)}A^7_{(110)}\la{)}\right.\\
&+\frac{8}{5}\la{(}A^4_{(020)}A^5_{(011)}-A^4_{(011)}A^5_{(020)}+
A^6_{(020)}A^7_{(011)}-A^6_{(011)}A^7_{(020)}\la{)}\\
&\left. -\frac{8}{5}\la{(}A^4_{(002)}A^5_{(011)}-A^4_{(011)}A^5_{(002)}+
A^6_{(002)}A^7_{(011)}-A^6_{(011)}A^7_{(002)}\la{)}\right]
\cdot \la{(} \int_{\mathbb{S}^2}\phi^2 \la{)}.\end{array}
\end{eqnarray*}}
In $(21)=(16)$ we have the following terms:
The terms with $a+b=(3,0,1)$ that are given by $a=(2,0,0)$ with $b=(1,0,1)$
and vice versa.
The terms with $a+b=(1,2,1)$ that are given by $a=(0,2,0)$ with
$b=(1,0,1)$, and  $a=(1,1,0)$ with $b=(0,1,1)$, and vice versa.
The terms with $a+b= (1,0,3)$ that are given by $a=(1,1,0)$ with $b=(0,1,1)$
and vice versa.
Thus,
{\small
\begin{eqnarray*}
\lefteqn{(16)=}\\
&&\begin{array}{l}
= 4 (A^4_{(200)}A^6_{(101)}-A^4_{(101)}A^6_{(200)}
-A^5_{(200)}A^7_{(101)}+A^5_{(101)}A^7_{(200)})\cdot\La{(}
-2\int_{\mathbb{S}^2}\phi_1^2\phi_2^2\La{)}\\
+  4 (A^4_{(110)}A^6_{(011)}-A^4_{(011)}A^6_{(110)}
-A^5_{(110)}A^7_{(011)}+A^5_{(011)}A^7_{(110)})
\cdot \La{(} \int_{\mathbb{S}^2}\phi_1^2\phi_2^2
+\int_{\mathbb{S}^2}\phi_2^2\phi_3^2- 
\int_{\mathbb{S}^2}\phi_2^4 \La{)}\\
+  4 (A^4_{(020)}A^6_{(101)}-A^4_{(101)}A^6_{(020)}
-A^5_{(020)}A^7_{(101)}+A^5_{(101)}A^7_{(020)})\cdot
\La{(} 2\int_{\mathbb{S}^2}\phi_1^2\phi_2^2-
2\int_{\mathbb{S}^2}\phi_2^2\phi_3^2\La{)}\\
+  4 (A^4_{(002)}A^6_{(101)}-A^4_{(101)}A^6_{(002)}
-A^5_{(002)}A^7_{(101)}+A^5_{(101)}A^7_{(002)})\cdot
 \La{(} 2\int_{\mathbb{S}^2}\phi_2^2\phi_3^2\La{)}\end{array}
\end{eqnarray*}}
that is, 
{\small
\begin{eqnarray*}
\begin{array}{ll}
(16) =&
\left[ -\frac{8}{5} \la{(}A^4_{(200)}A^6_{(101)}-A^4_{(101)}A^6_{(200)}
-A^5_{(200)}A^7_{(101)}+A^5_{(101)}A^7_{(200)}\la{)}\right.\\
& -\frac{4}{5} \la{(} A^4_{(110)}A^6_{(011)}-A^4_{(011)}A^6_{(110)}
-A^5_{(110)}A^7_{(011)}+A^5_{(011)}A^7_{(110)}\la{)}\\
&\left.  +\frac{8}{5} \la{(}A^4_{(002)}A^6_{(101)}-A^4_{(101)}A^6_{(002)}
-A^5_{(002)}A^7_{(101)}+A^5_{(101)}A^7_{(002)}\la{)}\right]
\cdot \la{(} \int_{\mathbb{S}^2}\phi^2\la{)}.
\end{array}
\end{eqnarray*}}
In $(22)=(17)$ we have the following terms:
The terms with $a+b=(3,1,0)$ that are given by $a=(2,0,0)$ with $b=(1,1,0)$
and vice versa.
The terms with $a+b=(1,3,0)$ that are given by $a=(0,2,0)$ with
$b=(1,1,0)$ and vice versa.
The terms with $a+b= (1,1,2)$ that are given by $a=(1,0,1)$ with $b=(0,1,1)$, 
and by $a=(0,0,2)$ with $(1,1,0)$, and vice versa.
Therefore, 
{\small
\begin{eqnarray*}
\lefteqn{(17)=}\\
&&\begin{array}{l}
=  4 (-A^5_{(200)}A^6_{(110)}+A^5_{(110)}A^6_{(200)}
-A^4_{(200)}A^7_{(110)}+A^4_{(110)}A^7_{(200)})\cdot
\La{(}2\int_{\mathbb{S}^2}\phi_1^2\phi_3^2\La{)}\\
+ 4 (-A^5_{(020)}A^6_{(110)}+A^5_{(110)}A^6_{(020)}
-A^4_{(020)}A^7_{(110)}+A^4_{(110)}A^7_{(020)})\cdot
\La{(} -2 \int_{\mathbb{S}^2}\phi_2^2\phi_3^2 \La{)}\\
4 (-A^5_{(101)}A^6_{(011)}+A^5_{(011)}A^6_{(101)}
-A^4_{(101)}A^7_{(011)}+A^4_{(011)}A^7_{(101)})\cdot
\La{(} -\int_{\mathbb{S}^2}\phi_2^2\phi_3^2
-\int_{\mathbb{S}^2}\phi_1^2\phi_3^2
+\int_{\mathbb{S}^2}\phi_3^4 \La{)}\\
  4 (-A^5_{(002)}A^6_{(110)}+A^5_{(110)}A^6_{(002)}
-A^4_{(002)}A^7_{(110)}+A^4_{(110)}A^7_{(002)})\cdot
\La{(} \int_{\mathbb{S}^2}\phi_2^2\phi_3^2
-\int_{\mathbb{S}^2}\phi_1^2\phi_3^2\La{)},\end{array}
\end{eqnarray*}}
that is
{\small
\[
\begin{array}{ll}
(17)=&\left[\frac{8}{5}\la{(}
-A^5_{(200)}A^6_{(110)}+A^5_{(110)}A^6_{(200)}
-A^4_{(200)}A^7_{(110)}+A^4_{(110)}A^7_{(200)}\la{)}\right.\\
&-\frac{8}{5}\la{(}
-A^5_{(020)}A^6_{(110)}+A^5_{(110)}A^6_{(020)}
-A^4_{(020)}A^7_{(110)}+A^4_{(110)}A^7_{(020)}\la{)}\\
&\left.+\frac{4}{5}\la{(}
-A^5_{(101)}A^6_{(011)}+A^5_{(011)}A^6_{(101)}
-A^4_{(101)}A^7_{(011)}+A^4_{(011)}A^7_{(101)}\la{)}\right]
\cdot \left(\int_{\mathbb{S}^2}\phi^2\right).
\end{array}
\] }
On the other hand, 
{\small
\begin{equation}
\begin{array}{ll}
\sum_{\alpha}\|\nabla f_{\alpha}\|^2 =& \left(\frac{8}{5} 
\sum_{\alpha}\la{(} (A^{\alpha}_{(200)})^2+(A^{\alpha}_{(020)})^2
+(A^{\alpha}_{(002)})^2 \La{)}\right. \\
&-\frac{4}{5}\sum_{\alpha} \la{(}A^{\alpha}_{(200)}A^{\alpha}_{(020)}
+A^{\alpha}_{(200)}A^{\alpha}_{(002)}+A^{\alpha}_{(020)}A^{\alpha}_{(002)}
\la{)}\\
&\left. +\frac{6}{5}\sum_{\alpha}\la{(} (A^{\alpha}_{(110)})^2+
 (A^{\alpha}_{(101)})^2+ (A^{\alpha}_{(011)})^2\La{)}\right)\cdot
\left(\int_{\mathbb{S}^2}\phi^2\right).
\end{array}\\[2mm]
\end{equation}
}
Thus, we have to show that $(15)+(16)+(17)\leq (24)$. This can be shown
 by proving that the following four inequalities hold
{\small 
$$\begin{array}{l}
4A^4_{(110)}A^5_{(101)}-8A^6_{(011)}A^7_{(020)}+8A^6_{(011)}A^7_{7(002)}
-8A^5_{(101)}A^7_{(200)}-4A^4_{(110)}A^6_{(011)}+8A^5_{(101)}A^7_{(002)}\\
+8A^4_{(110)}A^7_{(200)}-8A^4_{(110)}A^7_{(020)}-4A^5_{(101)}A^6_{(011)}
+4A^7_{(200)}A^7_{(020)}+4A^7_{(200)}A^7_{(002)}+4A^7_{(020)}A^7_{(002)}
\\[3mm]
\leq 8(A^7_{(200)})^2+8(A^7_{(020)})^2+8(A^7_{(002)})^2+
6(A^4_{4(110)})^2+6(A^5_{(101)})^2+6(A^6_{6(011)})^2,\\[2mm]
\end{array}
$$}
{\small
$$\begin{array}{l}
-4A^4_{(101)}A^5_{(110)}+8A^6_{(020)}A^7_{(011)}-8A^6_{(020)}A^7_{(011)}
+8A^4_{(101)}A^6_{(200)}+4A^5_{(110)}A^7_{(011)}-8A^4_{(101)}A^6_{(002)}\\
+8A^5_{(110)}A^6_{(200)}-8A^5_{(110)}A^6_{(020)}-4A^4_{(101)}A^7_{(011)}
+4A^6_{(200)}A^6_{(020)}+4A^6_{(200)}A^6_{(002)}
+4A^6_{(020)}A^6_{(002)}\\[3mm]
\leq 8(A^6_{6(200)})^2+8(A^6_{(020)})^2+8(A^6_{(002)})^2+
6(A^4_{(101)})^2+6(A^5_{(110)})^2 +6(A^7_{(011)})^2,\\[2mm]
\end{array}
$$}
{\small
$$\begin{array}{l}
4A^6_{(110)}A^7_{(101)}-8A^4_{(011)}A^5_{(020)}+8A^4_{(011)}A^5_{(002)}
+8A^5_{(200)}A^6_{(200)}+4A^4_{(011)}A^6_{(110)}-8A^5_{(002)}A^7_{(101)}\\
-8A^5_{(200)}A^6_{(110)}+8A^5_{(020)}A^6_{(110)}+4A^4_{(011)}A^7_{(101)}
+4A^5_{(200)}A^5_{(020)}+4A^5_{(200)}A^5_{(002)}+4A^5_{(020)}A^5_{(002)}
\\[3mm]
\leq 8(A^5_{(200)})^2+8(A^5_{(020)})^2+8(A^5_{(002)})^2+
6(A^4_{(011)})^2+6(A^6_{(110)})^2
+6(A^7_{(101)})^2,\\[2mm]
\end{array}
$$}
{\small
$$\begin{array}{l}
-4A^6_{(101)}A^7_{(110)}+8A^4_{(020)}A^5_{(011)}-8A^4_{(002)}A^5_{(011)}
-8A^4_{(200)}A^6_{(101)}-4A^5_{(011)}A^7_{(110)}+8A^4_{(002)}A^6_{(101)}\\
-8A^4_{(200)}A^7_{(110)}+8A^4_{(020)}A^7_{(110)}+4A^5_{(011)}A^6_{(101)}
+4A^4_{(200)}A^4_{(020)}+4A^4_{(200)}A^4_{(002)}+4A^4_{(020)}A^4_{(002)}
\\[3mm]
\leq 8(A^4_{(200)})^2+8(A^4_{(020)})^2+8(A^4_{(002)})^2+6(A^5_{(011)})^2
+6(A^6_{(101)})^2+6(A^7_{(110)})^2.
\end{array}\\[2mm]
$$}
The above inequalities are all of the form
\begin{equation}
\begin{array}{l}
2AB+2AC+2BC+2XY+2XZ+2YZ +4\langle(A,B,C),(X-Y,Y-Z,Z-X)\rangle \\
\quad \leq  3(A^2+B^2+C^2) +4(X^2+Y^2+Z^2)\end{array}
\end{equation}
Next lemma shows that (25) holds,  with equality iff $X=Y=Z=A=B=C=0$,
what proves the lemma. \qed
\begin{lemma} For any real numbers $A,B,C,X,Y,Z$, we have
\begin{equation}
\begin{array}{l}
2AB+2AC+2BC+2XY+2XZ+2YZ+4\langle (A,C,B),(X-Y,Y-Z,Z-X)\rangle\\[2mm]
\quad \quad \leq 3(A^2+B^2+C^2) +3(X^2+Y^2+Z^2), \end{array}
\end{equation}
with equality iff $(A,C,B)$ and $(X-Y,Y-Z,Z-X)$ are collinear.
\end{lemma}
\noindent
\em Proof. \em 
Note first that
$$\langle (A,C,B),(X-Y,Y-Z,Z-X)\rangle=\langle (X,Z,Y),(A-B,B-C,C-A)\rangle,$$
and $2|\langle u,v\rangle| \leq |u|^2+|v|^2$ with equality
iff $u,v$ are collinear. Then
\begin{eqnarray*}
\lefteqn{4\langle (A,C,B),(X-Y,Y-Z,Z-X)\rangle=}\\[1mm]
&=& 2 \langle (A,C,B),(X-Y,Y-Z,Z-X)\rangle \\
&& + 2 \langle (X,Z,Y),(A-B,B-C,C-A)\rangle \\[1mm]
&\leq & A^2+B^2+C^2+ (X-Y)^2+(Y-Z)^2 +(Z-X)^2\\
&& + X^2+Y^2+Z^2 +(A-B)^2+(B-C)^2+(C-A)^2.
\end{eqnarray*}
But this is just (26).\qed
\begin{lemma} If $l=3$ the functions
$$\begin{array}{lll}
f_4= 3\phi_3^3+ 4\phi_1^2\phi_3 +4\phi_2^2\phi_3 & & f_5=
-4\phi_2^3-4\phi_2\phi_3^2-3\phi_1^2\phi_2\\
f_6=-\phi_1^3-\phi_1\phi_2^2-2\phi_1\phi_3^2& & f_7=\phi_1\phi_2\phi_3
\end{array}$$
do not satisfy the inequality (3).
\end{lemma}
\noindent
\em Proof. \em
We will show how we have found such functions. We have
$$\begin{array}{l}
\int_{\mathbb{S}^2}\phi^{(6,0,0)}dM=
\frac{3}{7}\int_{\mathbb{S}^2}\phi^2dM \\[1mm]
\int_{\mathbb{S}^2}\phi^{(4,2,0)}dM=
\frac{3}{5\times 7}\int_{\mathbb{S}^2}\phi^2dM \\[1mm]
\int_{\mathbb{S}^2}\phi^{(2,2,2)}dM=
\frac{1}{5\times 7}\int_{\mathbb{S}^2}\phi^2dM.
\end{array}$$
The multi-powers $a,b$ are running all over 
$$\{(300),(030),(003),(210),(201),(021),(120),
(012),(102),(111)\}.$$
The terms with $a+b=b+a=(E,O,O)$ are given by $(1,1,1)+(3,0,0)$, 
$(0,0,3)+(0,3,0)$, $(2,0,1)+(0,3,0)$, $(0,2,1)+(0,3,0)$,
$(2,1,0)+(0,0,3)$, $(0,1,2)+(0,0,3)$, $(2,0,1)+(2,1,0)$, $(0,2,1)+(2,1,0)$,
$(0,1,2)+(2,0,1)$, $(1,1,1)+(1,2,0)$, $(1,1,1)+(1,0,2)$. Thus,
{\small
\begin{eqnarray*}
\lefteqn{(15)=}\\
&&\!\!\!\!\!\!\!\!\!\!\!\!\!\!\begin{array}{l}
= 4\la{(}A^4_{(111)}A^5_{(300)}-A^4_{(300)}A^5_{(111)}
+A^6_{(111)}A^7_{(300)}-A^6_{(300)}A^7_{(111)}\la{)}\cdot
\left(-3\int_{\mathbb{S}^2}\phi_1^4\phi_3^2 + 
3\int_{\mathbb{S}^2}\phi_1^4\phi_2^2 \right)\\
+  4\la{(}A^4_{(003)}A^5_{(030)}-A^4_{(030)}A^5_{(003)}
+A^6_{(003)}A^7_{(030)}-A^6_{(030)}A^7_{(003)}\la{)}
\cdot \left(-9\int_{\mathbb{S}^2}\phi_1^2\phi_2^2\phi_3^2 \right)\\
+  4\la{(}A^4_{(201)}A^5_{(030)}-A^4_{(030)}A^5_{(201)}
+A^6_{(201)}A^7_{(030)}-A^6_{(030)}A^7_{(201)}\la{)}\cdot
\left(6\int_{\mathbb{S}^2}\phi_1^2\phi_2^2\phi_3^2
-3\int_{\mathbb{S}^2}\phi_1^4\phi_2^2 \right)\\
+  4\la{(}A^4_{(021)}A^5_{(030)}-A^4_{(030)}A^5_{(021)}
+A^6_{(021)}A^7_{(030)}-A^6_{(030)}A^7_{(021)}\la{)}\cdot
\left(-3\int_{\mathbb{S}^2}\phi_1^2\phi_2^4 \right)\\
+  4\la{(}A^4_{(012)}A^5_{(003)}-A^4_{(003)}A^5_{(012)}
+A^6_{(012)}A^7_{(003)}-A^6_{(003)}A^7_{(120)}\la{)}
\cdot \left(3\int_{\mathbb{S}^2}\phi_1^2\phi_3^4 \right)\\
+  4\la{(}A^4_{(201)}A^5_{(210)}-A^4_{(210)}A^5_{(201)}
+A^6_{(201)}A^7_{(210)}-A^6_{(210)}A^7_{(201)}\la{)}\cdot
\left(2\int_{\mathbb{S}^2}\phi_1^4\phi_3^2 + 
2\int_{\mathbb{S}^2}\phi_1^4\phi_2^2
-\int_{\mathbb{S}^2}\phi_1^6 \right)\\
+  4\la{(}A^4_{(021)}A^5_{(210)}-A^4_{(210)}A^5_{(021)}
+A^6_{(021)}A^7_{(210)}-A^6_{(210)}A^7_{(021)}\la{)}\cdot
\left(-4\int_{\mathbb{S}^2}\phi_1^2\phi_2^2\phi_3^2 + 
2\int_{\mathbb{S}^2}\phi_1^2\phi_2^4
-\int_{\mathbb{S}^2}\phi_1^4\phi_2^2  \right)\\
+  4\la{(}A^4_{(012)}A^5_{(201)}-A^4_{(201)}A^5_{(012)}
+A^6_{(012)}A^7_{(201)}-A^6_{(201)}A^7_{(012)}\la{)}\cdot
\left(-2\int_{\mathbb{S}^2}\phi_1^2\phi_3^4 + 
4\int_{\mathbb{S}^2}\phi_1^2\phi_2^2\phi_3^2 
+\int_{\mathbb{S}^2}\phi_1^4\phi_3^2  \right)\\
+  4\la{(}A^4_{(111)}A^5_{(120)}-A^4_{(120)}A^5_{(111)}
+A^6_{(111)}A^7_{(120)}-A^6_{(120)}A^7_{(111)}\la{)}
\cdot \left(\int_{\mathbb{S}^2}\phi_1^2\phi_2^2\phi_3^2 + 
\int_{\mathbb{S}^2}\phi_1^2\phi_2^4 
-2\int_{\mathbb{S}^2}\phi_1^4\phi_2^2  \right)\\
+  4\la{(}A^4_{(111)}A^5_{(102)}-A^4_{(102)}A^5_{(111)}
+A^6_{(111)}A^7_{(102)}-A^6_{(102)}A^7_{(111)}\la{)}\cdot
\left(-\int_{\mathbb{S}^2}\phi_1^2\phi_3^4 
-\int_{\mathbb{S}^2}\phi_1^2\phi_2^2\phi_3^2 
+2\int_{\mathbb{S}^2}\phi_1^4\phi_3^2  \right)
\end{array}
\end{eqnarray*}}
The terms with $a+b=b+a=(O,E,O)$ are given by $(0,0,3)+(3,0,0)$, 
$(2,0,1)+(3,0,0)$, $(0,2,1)+(3,0,0)$, $(1,1,1)+(0,3,0)$,
$(1,2,0)+(0,0,3)$, $(1,0,2)+(0,0,3)$, $(1,0,2)+(2,1,0)$, $(1,1,1)+(2,1,0)$,
$(1,2,0)+(2,0,1)$, $(1,0,2)+(2,0,1$, $(1,2,0)+(0,2,1)$, $(1,0,2)+(0,2,1)$,
$(1,1,1)+(0,1,2)$. Hence, 
{\small 
\begin{eqnarray*}
\lefteqn{(16)=}\\
&&\!\!\!\!\!\!\!\!\!\!\!\!\!\!\!\begin{array}{l}
= 4\la{(}A^4_{(003)}A^6_{(300)}-A^4_{(300)}A^6_{(003)}
-A^5_{(003)}A^7_{(300)}+A^5_{(300)}A^7_{(003)}\la{)}
\cdot \left(9\int_{\mathbb{S}^2}\phi_1^2\phi_2^2\phi_3^2 \right)\\
+4\la{(}A^4_{(201)}A^6_{(300)}-A^4_{(300)}A^6_{(201)}
-A^5_{(201)}A^7_{(300)}+A^5_{(300)}A^7_{(201)}\la{)}
\cdot \left(3\int_{\mathbb{S}^2}\phi_1^4\phi_2^2dM \right)\\
+ 4\la{(}A^4_{(021)}A^6_{(300)}-A^4_{(300)}A^6_{(021)}
-A^5_{(021)}A^7_{(300)}+A^5_{(300)}A^7_{(021)}\la{)}\cdot
\left(-6\int_{\mathbb{S}^2}\phi_1^2\phi_2^2\phi_3^2
+3\int_{\mathbb{S}^2}\phi_1^2\phi_2^4 \right)\\
+ 4\la{(}A^4_{(111)}A^6_{(030)}-A^4_{(030)}A^6_{(111)}
-A^5_{(111)}A^7_{(030)}+A^5_{(030)}A^7_{(111)}\la{)}\cdot
\left(-3\int_{\mathbb{S}^2}\phi_1^2\phi_2^4 +
3\int_{\mathbb{S}^2}\phi_2^4\phi_3^2 \right)\\
+ 4\la{(}A^4_{(120)}A^6_{(003)}-A^4_{(003)}A^6_{(120)}
-A^5_{(120)}A^7_{(003)}+A^5_{(003)}A^7_{(120)}\la{)}
 \cdot \left(6\int_{\mathbb{S}^2}\phi_1^2\phi_2^2\phi_3^2
-3\int_{\mathbb{S}^2}\phi_2^4\phi_3^2 \right)
\end{array}
\end{eqnarray*}
}
{\small \begin{eqnarray*}\begin{array}{l}
+  4\la{(}A^4_{(111)}A^6_{(210)}-A^4_{(210)}A^6_{(111)}
-A^5_{(111)}A^7_{(210)}+A^5_{(210)}A^7_{(111)}\la{)}
 \cdot \left(-\int_{\mathbb{S}^2}\phi_1^4\phi_2^2
-\int_{\mathbb{S}^2}\phi_1^2\phi_2^2\phi_3^2
+2\int_{\mathbb{S}^2}\phi_1^2\phi_2^4 \right)\\
+  4\la{(}A^4_{(120)}A^6_{(201)}-A^4_{(201)}A^6_{(120)}
-A^5_{(120)}A^7_{(201)}+A^5_{(201)}A^7_{(120)}\la{)}
\cdot \left(2\int_{\mathbb{S}^2}\phi_1^4\phi_2^2
-4\int_{\mathbb{S}^2}\phi_1^2\phi_2^2\phi_3^2 
-\int_{\mathbb{S}^2}\phi_1^2\phi_2^4  \right)\\
+  4\la{(}A^4_{(102)}A^6_{(201)}-A^4_{(201)}A^6_{(102)}
-A^5_{(102)}A^7_{(201)}+A^5_{(201)}A^7_{(102)}\la{)}\cdot
\left(3\int_{\mathbb{S}^2}\phi_1^2\phi_2^2\phi_3^2  \right)\\
+  4\la{(}A^4_{(120)}A^6_{(021)}-A^4_{(021)}A^6_{(120)}
-A^5_{(120)}A^7_{(021)}+A^5_{(021)}A^7_{(120)}\la{)}
\cdot \left(2\int_{\mathbb{S}^2}\phi_1^2\phi_2^4
+2\int_{\mathbb{S}^2}\phi_2^4\phi_3^2
-\int_{\mathbb{S}^2}\phi_2^6  \right)\\
+ 4\la{(}A^4_{(102)}A^6_{(021)}-A^4_{(021)}A^6_{(102)}
-A^5_{(102)}A^7_{(021)}+A^5_{(021)}A^7_{(102)}\la{)}
\cdot \left(-4\int_{\mathbb{S}^2}\phi_1^2\phi_2^2\phi_3^2  
+2\int_{\mathbb{S}^2}\phi_2^2\phi_3^4
-\int_{\mathbb{S}^2}\phi_2^4\phi_3^2  \right)\\
+  4\la{(}A^4_{(111)}A^6_{(012)}-A^4_{(012)}A^6_{(111)}
-A^5_{(111)}A^7_{(012)}+A^5_{(012)}A^7_{(111)}\la{)}\cdot
\left(\int_{\mathbb{S}^2}\phi_1^2\phi^2\phi_3^2 
+\int_{\mathbb{S}^2}\phi_2^2\phi_3^4 
-2\int_{\mathbb{S}^2}\phi_2^4\phi_3^2  \right)
\end{array}
\end{eqnarray*}}

The terms with $a+b=b+a= (O,O,E)$ are given by $(0,3,0)+(3,0,0)$
$(2,1,0)+(3,0,0)$, $(0,1,2)+(3,0,0)$, $(1,2,0)+(0,3,0)$,
$(1,0,2)+(0,3,0)$, $(1,1,1)+(0,0,3)$, $(1,2,0)+(2,1,0)$, $(1,1,1)+(2,0,1)$,
$(1,1,1)+(0,2,1)$, $(0,1,2)+(1,2,0)$, $(1,0,2)+(0,1,2)$. Therefore
{\small
\begin{eqnarray*}
\lefteqn{(17)=}\\
&&\!\!\!\!\!\!\!\!\!\!\!\!\!\!\!\!\!\begin{array}{l}
= 4\la{(}-A^4_{(030)}A^7_{(300)}+A^4_{(300)}A^7_{(003)}
-A^5_{(030)}A^6_{(300)}+A^5_{(300)}A^6_{(003)}\la{)}
\cdot \left(-9\int_{\mathbb{S}^2}\phi_1^2\phi_2^2\phi_3^2 \right)\\
+4\la{(}-A^4_{(210)}A^7_{(300)}+A^4_{(300)}A^7_{(210)}
-A^5_{(210)}A^6_{(300)}+A^5_{(300)}A^6_{(210)}\la{)}
\cdot \left(-3\int_{\mathbb{S}^2}\phi_1^4\phi_3^2 \right)\\
+ 4\la{(}-A^4_{(012)}A^7_{(300)}+A^4_{(300)}A^7_{(012)}
-A^5_{(012)}A^6_{(300)}+A^5_{(300)}A^6_{(012)}\la{)}\cdot
\left(6\int_{\mathbb{S}^2}\phi_1^2\phi_2^2\phi_3^2
-3\int_{\mathbb{S}^2}\phi_1^2\phi_3^4 \right)\\
+ 4\la{(}-A^4_{(120)}A^7_{(030)}+A^4_{(030)}A^7_{(120)}
-A^5_{(120)}A^6_{(030)}+A^5_{(030)}A^6_{(120)}\la{)}\cdot
\left(3\int_{\mathbb{S}^2}\phi_2^4\phi_3^2 \right)\\
+  4\la{(}-A^4_{(102)}A^7_{(030)}+A^4_{(030)}A^7_{(102)}
-A^5_{(102)}A^6_{(030)}+A^5_{(030)}A^6_{(102)}\la{)}\cdot
\left(-6\int_{\mathbb{S}^2}\phi_1^2\phi_2^2\phi_3^2
+ 3\int_{\mathbb{S}^2}\phi_2^2\phi_3^4 \right)\\
+  4\la{(}-A^4_{(111)}A^7_{(003)}+A^4_{(003)}A^7_{(111)}
-A^5_{(111)}A^6_{(003)}+A^5_{(003)}A^6_{(111)}\la{)}\cdot
\left(-3\int_{\mathbb{S}^2}\phi_2^2\phi_3^4 + 
3\int_{\mathbb{S}^2}\phi_1^2\phi_3^4   \right)\\
+  4\la{(}-A^4_{(120)}A^7_{(210)}+A^4_{(210)}A^7_{(120)}
-A^5_{(120)}A^6_{(210)}+A^5_{(210)}A^6_{(120)}\la{)}\cdot
\left(-3\int_{\mathbb{S}^2}\phi_1^2\phi_2^2\phi_3^2 \right)\\
+  4\la{(}-A^4_{(111)}A^7_{(201)}+A^4_{(201)}A^7_{(111)}
-A^5_{(111)}A^6_{(201)}+A^5_{(201)}A^6_{(111)}\la{)}\cdot
 \left(\int_{\mathbb{S}^2}\phi_1^2\phi_2^2\phi_3^2  + 
\int_{\mathbb{S}^2}\phi_1^4\phi_3^2
-2\int_{\mathbb{S}^2}\phi_1^2\phi_3^4  \right)\\
+  4\la{(}-A^4_{(111)}A^7_{(021)}+A^4_{(021)}A^7_{(111)}
-A^5_{(111)}A^6_{(021)}+A^5_{(021)}A^6_{(111)}\la{)}
\cdot \left(-\int_{\mathbb{S}^2}\phi_2^4\phi_3^2
-\int_{\mathbb{S}^2}\phi_1^2\phi_2^2\phi_3^2  
+2\int_{\mathbb{S}^2}\phi_2^2\phi_3^4  \right)\\
+  4\la{(}-A^4_{(012)}A^7_{(120)}+A^4_{(120)}A^7_{(012)}
-A^5_{(012)}A^6_{(120)}+A^5_{(120)}A^6_{(012)}\la{)}\cdot
\left(2\int_{\mathbb{S}^2}\phi_2^4\phi_3^2  
-4\int_{\mathbb{S}^2}\phi_1^2\phi_2^2\phi_3^2
-\int_{\mathbb{S}^2}\phi_2^2\phi_3^4  \right)\\
+  4\la{(}-A^4_{(102)}A^7_{(210)}+A^4_{(210)}A^7_{(102)}
-A^5_{(102)}A^6_{(210)}+A^5_{(210)}A^6_{(102)}\la{)}\cdot
\left(4\int_{\mathbb{S}^2}\phi_1^2\phi_2^2\phi_3^2  + 
-2\int_{\mathbb{S}^2}\phi_1^4\phi_3^2 
+ \int_{\mathbb{S}^2}\phi_1^2\phi_3^4  \right)
\end{array}
\end{eqnarray*}}
On the other hand
{\small
\begin{eqnarray}
\lefteqn{\sum_{\alpha}\|\nabla f_{\alpha}\|^2_{L^2}=} \nonumber\\[-3mm]
&&\begin{array}{l}
(\int_{\mathbb{S}^2}\phi^2dM)
\cdot \left(\sum_{\alpha}\frac{9\times 6}{5\times 7}
\la{(} (A^{\alpha}_{(300)})^2
+ (A^{\alpha}_{(030)})^2+ (A^{\alpha}_{(003)})^2\la{)}\right. \\
+\sum_{\alpha}\frac{11\times 2}{5\times 7}\la{(}(A^{\alpha}_{(210)})^2+
(A^{\alpha}_{(201)})^2+(A^{\alpha}_{(021)})^2+(A^{\alpha}_{(120)})^2
+(A^{\alpha}_{(012)})^2+(A^{\alpha}_{(1,0,2)})^2\la{)}\\
 +\sum_{\alpha}\frac{3\times 4}{5\times 7}(A^{\alpha}_{(111)})^2
-\sum_{\alpha}\frac{ 4}{5\times 7}\la{(}
A^{\alpha}_{((210)}A^{\alpha}_{(012)}
+A^{\alpha}_{(201)}A^{\alpha}_{(021)}
+A^{\alpha}_{(120)}A^{\alpha}_{(102)}\la{)}\\
\left. -\sum_{\alpha}\frac{12}{5\times 7}\la{(}
 A^{\alpha}_{(300)}A^{\alpha}_{(120)}
+A^{\alpha}_{(300)}A^{\alpha}_{(102)}+A^{\alpha}_{(030)}A^{\alpha}_{(210)}
+A^{\alpha}_{(030)}A^{\alpha}_{(012)}+
A^{\alpha}_{(003)}A^{\alpha}_{(201)}
+A^{\alpha}_{(003)}A^{\alpha}_{(021)}\la{)}\right)
\end{array}
\end{eqnarray}
}
Then $(15)+(16)+(17)\leq (27)$ to hold for all possible coefficients
$A^{\alpha}_a$ is equivalent
to  two linearly  independent inequalities:
\newpage
{\small
\begin{eqnarray}
\begin{array}{l}
-36A^4_{(003)}A^5_{(030)}-12A^4_{(201)}A^5_{(030)}-36A^4_{(021)}A^5_{(030)}
-36A^4_{(003)}A^5_{(012)}-12A^4_{(201)}A^5_{(210)}-4A^4_{(021)}A^5_{(210)}\\
-4A^4_{(201)}A^5_{(012)}+8A^6_{(120)}A^7_{(111)}-8A^6_{(102)}A^7_{(111)}
+36A^4_{(003)}A^6_{(300)}+36A^4_{(201)}A^6_{(300)}+12A^4_{(021)}A^6_{(300)}\\
+12A^4_{(003)}A^6_{(120)}+8A^5_{(210)}A^7_{(111)}+4A^4_{(201)}A^6_{(120)}
-12A^4_{(201)}A^6_{(102)}+12A^4_{(021)}A^6_{(120)}+4A^4_{(021)}A^6_{(102)}\\
-8A^5_{(012)}A^7_{(111)} + 36A^5_{(030)}A^6_{(300)}+ 36A^5_{(210)}A^6_{(300)}
+12A^5_{(012)}A^6_{(300)}+36A^5_{(030)}A^6_{(120)}+12A^5_{(030)}A^6_{(102)}\\
-12A^5_{(210)}A^6_{(120)}-8A^4_{(201)}A^7_{(111)}+8A^4_{(021)}A^7_{(111)}
+4A^5_{(012)}A^6_{(120)}+4A^5_{(210)}A^6_{(102)}\\
+ 4\la{(}A^4_{(201)}A^4_{(021)}+A^5_{(210)}A^5_{(012)}
+A^6_{(120)}A^6_{(102)}\la{)}\\
+12\la{(}A^6_{(300)}A^6_{(120)}+A^6_{(300)}A^6_{(102)}+A^5_{(030)}A^5_{(210)}
+A^5_{(030)}A^5_{(012)}+A^4_{(003)}A^4_{(201)}+ A^4_{(003)}A^4_{(021)}\la{)}
\\
\leq 
~~54\la{(} (A^4_{(003)})^2+(A^5_{(030)})^2+(A^6_{(300)})^2\la{)}
+12(A^7_{(111)})^2\\
~~~~+22\la{(} (A^4_{(201)})^2+(A^4_{(021)})^2+(A^5_{(012)})^2+(A^5_{(210)})^2
+(A^6_{(120)})^2+(A^6_{(102)})^2\la{)}
\end{array}
\end{eqnarray}
}
and 
{\small
\begin{eqnarray}
\begin{array}{l}
36A^4_{(030)}A^5_{(003)}-36A^6_{(003)}A^7_{(030)}+36A^6_{(030)}A^7_{(003)}
+12A^4_{(030)}A^5_{(201)}-12A^6_{(201)}A^7_{(030)}+12A^6_{(030)}A^7_{(201)}\\
+36A^4_{(030)}A^5_{(021)}-36A^6_{(021)}A^7_{(030)}+36A^6_{(030)}A^7_{(021)}
+36A^4_{(012)}A^5_{(003)}+36A^6_{(012)}A^7_{(003)}-36A^6_{(003)}A^7_{(120)}\\
+12A^4_{(210)}A^5_{(201)}-12A^6_{(201)}A^7_{(210)}+12A^6_{(210)}A^7_{(201)}
+4A^4_{(210)}A^5_{(021)}-4A^6_{(021)}A^7_{(210)}+4A^6_{(210)}A^7_{(021)}\\
+4A^4_{(012)}A^5_{(201)} + 4A^6_{(012)}A^7_{(201)}-4A^6_{(201)}A^7_{(012)}
-8A^4_{(111)}A^5_{(120)}+8A^4_{(120)}A^5_{(111)}-8A^6_{(111)}A^7_{(120)}\\
+8A^4_{(111)}A^5_{(102)}-8A^4_{(102)}A^5_{(111)}+8A^6_{(111)}A^7_{(102)}
-36A^4_{(300)}A^6_{(003)}-36A^5_{(003)}A^7_{(300)}+36A^6_{(300)}A^7_{(003)}\\
-36A^4_{(300)}A^6_{(201)}-36A^5_{(201)}A^7_{(300)}+36A^5_{(300)}A^7_{(201)}
-12A^4_{(300)}A^6_{(021)}-12A^5_{(021)}A^7_{(300)}+12A^5_{(300)}A^7_{(021)}\\
-12A^4_{(120)}A^6_{(003)}+12A^5_{(120)}A^7_{(003)}-12A^5_{(003)}A^7_{(120)}
+8A^4_{(111)}A^6_{(210)}-8A^4_{(210)}A^6_{(111)}-8A^5_{(111)}A^7_{(210)}\\
-4A^4_{(120)}A^6_{(201)}+4A^5_{(120)}A^7_{(201)}-4A^5_{(201)}A^7_{(120)}
+12A^4_{(102)}A^6_{(201)} - 12A^5_{(102)}A^7_{(201)}+12A^5_{(201)}A^7_{(102)}\\
-12A^4_{(120)}A^6_{(021)}+12A^5_{(120)}A^7_{(021)}-12A^5_{(021)}A^7_{(120)}
-4A^4_{(102)}A^6_{(021)}+4A^5_{(102)}A^7_{(021)}-4A^5_{(021)}A^7_{(102)}\\
-8A^4_{(111)}A^6_{(012)}+8A^4_{(012)}A^6_{(111)}+8A^5_{(111)}A^7_{(012)}
+36A^4_{(030)}A^7_{(300)}-36A^4_{(300)}A^7_{(003)}-36A^5_{(300)}A^6_{(003)}\\
+36A^4_{(210)}A^7_{(300)}-36A^4_{(300)}A^7_{(210)}-36A^5_{(300)}A^6_{(210)}
+12A^4_{(012)}A^7_{(300)}-12A^4_{(300)}A^7_{(012)}-12A^5_{(300)}A^6_{(012)}\\
-36A^4_{(120)}A^7_{(030)}+36A^4_{(030)}A^7_{(120)}-36A^5_{(120)}A^6_{(030)}
-12A^4_{(102)}A^7_{(030)}+12A^4_{(030)}A^7_{(102)}-12A^5_{(102)}A^6_{(030)}\\
+12A^4_{(120)}A^7_{(210)}-12A^4_{(210)}A^7_{(120)}+12A^5_{(120)}A^6_{(201)}
+8A^4_{(111)}A^7_{(201)}+8A^5_{(111)}A^6_{(201)}-8A^5_{(201)}A^6_{(111)}\\
-8A^4_{(111)}A^7_{(021)}-8A^5_{(111)}A^6_{(021)}+8A^5_{(021)}A^6_{(111)}
+4A^4_{(012)}A^7_{(120)}-4A^4_{(120)}A^7_{(012)}-4A^5_{(120)}A^6_{(012)}\\
-4A^4_{(102)}A^7_{(201)}+4A^5_{(210)}A^7_{(102)}-4A^5_{(102)}A^6_{(210)}\\
+4\la{(}
+A^4_{(210)}A^4_{(012)}+A^6_{(210)}A^6_{(012)}+A^7_{(210)}A^7_{(012)}
+A^5_{(201)}A^5_{(021)}+ A^6_{(201)}A^6_{(021)}\la{)}\\
+4\la{(} A^7_{(201)}A^7_{(021)}
+A^4_{(120)}A^4_{(102)}+ A^5_{(120)}A^5_{(102)}+A^7_{(120)}A^7_{(102)}\la{)}\\
+12\la{(}A^4_{(300)}A^4_{(120)}+A^4_{(300)}A^4_{(102)}+
A^5_{(300)}A^5_{(120)}+A^5_{(300)}A^5_{(102)}
+A^7_{(300)}A^7_{(120)}+A^7_{(300)}A^7_{(102)}\la{)}\\
+12\la{(}A^4_{(030)}A^4_{(210)}+A^4_{(030)}A^4_{(012)}
+A^6_{(030)}A^6_{(210)}+A^6_{(030)}A^6_{(012)}
+A^7_{(030)}A^7_{(210)}+A^7_{(030)}A^7_{(012)}\la{)}\\
+12\la{(}+A^5_{(003)}A^5_{(201)}+ A^5_{(003)}A^5_{(021)}
+A^6_{(003)}A^6_{(201)}+ A^6_{(003)}A^6{(021)}
+A^7_{(003)}A^7_{(201)}+ A^7_{(003)}A^7_{(021)}\la{)}\\
\leq 54\la{(} (A^4_{(300)})^2+  (A^4_{(030)})^2+ (A^5_{(300)})^2
+(A^5_{(003)})^2+(A^6_{(030)})^2+\\
~~~~~~~~~ +(A^6_{(003)})^2
+(A^7_{(003)})^2+(A^7_{(030)})^2 +(A^7_{(003)})^2
\la{)}
~~~+12\la{(}(A^4_{(111)})^2+(A^5_{(111)})^2+(A^6_{(111)})^2\la{)}\\
~~~~+22\la{(} 
(A^4_{(210)})^2+(A^4_{(120)})^2+(A^4_{(102)})^2+(A^4_{(012)})^2
+A^5_{(201)})^2+(A^5_{(120)})^2+(A^5_{(102)})^2+(A^5_{(021)})^2\la{)}\\
~~~~+\la{(}(A^6_{(210)})^2+(A^6_{(201)})^2+ (A^6_{(021)})^2+(A^6_{(012)})^2
\la{)}\\
~~~+\la{(}(A^7_{(210)})^2+(A^7_{(201)})^2
(A^7_{(120)})^2+(A^7_{(102)})^2+(A^7_{(021)})^2+(A^7_{(012)})^2
\la{)}
\end{array}
\end{eqnarray}
}
Now we consider only (28).
We define the polynomial function on 10 variables

$$\begin{array}{l}
F[a, b, c, x, y, z, w, u, v, g] := \\[2mm]
~~~~~ 54 (a^2 + b^2 + c^2) + 22 (x^2 + y^2 + z^2 + w^2 + u^2 + v^2) + 12g^2 \\
~~~~~  + 36 a b - 36 a c - 36 b c + 36 y b + 36 a z - 36 x c \\ 
~~~~~  -36 w c - 36 b u + 12 x b + 12 x w - 12 y c - 12 a u + 12 x v \\ 
~~~~~  -12 y u - 12 z c - 12 b v + 12 w u - 12 c u - 12 c v - 12 b w \\
~~~~~  -12 b z - 12 a x - 12 a y - 8 g (u - v + w - x + y - z)\\ 
~~~~~  +4 y w +4 x z - 4 x u - 4 y v - 4 z u - 4 w v - 4 x y - 4 w z - 4 u v
\end{array}$$
The inequality (28) is equivalent to say $F[a,b,c,x,y,w,u,v]\geq 0$
where 
$$\begin{array}{lllllll}
a=A^4_{(003)} & b= A^5_{(030)}& c=A^6_{(300)}& g=A^7_{(111)}\\
x=A^4_{(201)} & y= A^4_{(021)}& z=A^5_{(012)}& w=A^5_{(210)}
&u=A^6_{(120)}& v=A^6_{(102)}\end{array}$$
Using the Mathematica programming
we see that $F$ has only a critical point, that is at $0$, 
where $F$ vanish.
Then we compute
the Hessian of $F$, giving a $10\times 10$ matrix. Using a 
"Diag" package in Fortran 77 programming ( see \cite{Hahn}), we 
are able to obtain the eigenvalues of $\mathrm{Hess} F$ at the
origin, obtaining
$$\begin{array}{l}
193.95260118883090 \\       
111.22289635621148 \\       
123.94135950568288 \\      
-3.6844648605223074 \\       
64.826325872315152  \\      
39.493408799262383  \\      
5.8125282188336085  \\      
26.731868670430774  \\      
34.522364101735334  \\      
15.181112147219768 \\
\end{array}$$
Then we see there is a negative eigenvalue, what shows that there is
a direction given by the eigenvector corresponding to this negative
eigenvalue where $F$ might be negative.
Using the same Fortran programming we obtain the eigenvectors
{\tiny
$$\begin{array}{l}
(0.50859017, 0.56532378, -0.57392096, 0.15439833, 0.15001230, 0.12873210, 
9.32173475E-02, -0.14336746,     1.08674619E-02, 9.13834592E-03)\\[2mm]
(-0.70249306,     0.29659516,    -0.22777210,    0.28225840,     0.32578962,    -0.38748263,     
0.16652590 ,   -8.72224667E-03, 2.92987380E-03, -5.37370482E-02 )\\[2mm]
(2.63951580E-02, 0.52362033,     0.53104930,    -0.24485076,   0.21059064,    -0.11992104,
    -0.38731071,    -0.38285444,    -0.16684786,     2.23770015E-03)\\[2mm]
( 0.35473784,    -0.43567216,    -4.01160688E-02, 0.34413141,     0.44414138,    -0.46720073,
    -0.28168146,    -7.49844407E-02, -0.21986490,     0.12437580   )\\[2mm] 
( 0.12019036,     2.41075175E-02,-2.91444063E-02, -0.44169635,     0.19279750,    -0.15999765,
     0.32263696,     0.36108327,    -0.57445451,    -0.40222405    )\\[2mm]
(-0.17223047,    -0.21619180,     1.17000685E-02, -9.76854658E-03, 0.61789725,     0.69109876,
     9.10242003E-02, -0.20175159,    -0.11895024,      2.85076642E-02)\\[2mm]
( 0.20780841,     0.23610482 ,    0.50748260,     0.33018158,     0.23542015,     4.49400444E-02,
 0.33261644,     0.52470301,     0.19633343,    0.22929391  )\\[2mm]  
(-0.14729467,     0.13044269,    -0.22190370,    -7.04617534E-02, 8.01743129E-02, 0.18808081,  
  -0.68512740,     0.61403882,    -6.77787742E-02, 0.11934559)\\[2mm]    
( 0.10034294,    -9.04353313E-02, -7.49166773E-02, -0.44878484,     0.38166992,    -0.19261905 ,
   -3.88942769E-02, 5.71153059E-02, 0.72842332,    -0.23787165 )\\[2mm]   
(-4.73928140E-02, -1.38729246E-02, -0.14742221,     -0.45555392,     5.70707731E-02,
-0.15137651,     0.20074756,    -2.15775514E-02, -7.61975838E-02, 0.83399977 )   
\end{array}$$
}
where $E-0k$ means that we have to multiply the number by $10^{-k}$.
Then the eigenvector corresponding to the negative eigenvalue is
$ S=[a,b,c,x,y,z,w,u,v,g]$ where
$$\begin{array}{lllll} 
a=0.35473784 & & b= - 0.43567216 & & c= - 0.0401160688\\
 x=   0.34413141 & & y=   0.44414138 & & z=- 0.46720073\\
w= - 0.28168146 & & u= - 0.0749844407 & &  v=- 0.21986490\\
g=  0.12437580 \end{array} $$
We take the vector $[3, -4, -1, 3, 4, -4, -3, -1, -2, 1]$ in a neigbourhood of  $10 S$, and verify that
$F[3, -4, -1, 3, 4, -4, -3, -1, -2, 1]=-138 <0,  $
what proves our lemma.\qed

\section{Acknowledgments}
The author thanks George Rupp for computational assistance 
on Fortran programming.
 
\end{document}